\documentclass{article}

\usepackage{todonotes}
\usepackage{geometry}
\usepackage{amsmath}
\usepackage{amssymb}
\usepackage{amsthm}
\usepackage{bm}
\usepackage{graphicx}
\usepackage{color}
\usepackage{tabu}
\usepackage{booktabs}
\usepackage{caption}
\usepackage{authblk}
\usepackage{cite}
\usepackage{subfigure}
\usepackage{enumerate}
\usepackage{tabu}
\usepackage{multirow}
\usepackage{diagbox}
\usepackage[ruled, linesnumbered]{algorithm2e}

\newtheorem{defi}{Definition}

\newcommand{\eps}{\varepsilon}
\newcommand{\R}{\mathbb{R}}

\newcommand{\sS}{\mathbb{S}}
\newcommand{\bn}{\bm{n}}
\newcommand{\bx}{\bm{x}}
\newcommand{\by}{\bm{y}}

\newcommand{\bw}{\bm{w}}
\newcommand{\bv}{\bm{v}}
\newcommand{\bmu}{\bm{\mu}}
\newcommand{\bsigma}{\bm{\sigma}}
\newcommand{\bA}{\bm{A}}
\newcommand{\LO}{\mathcal{L}}
\newcommand*{\diff}{\mathop{}\!\mathrm{d}}

\newcommand{\avg}[1]{\left \langle #1 \right \rangle}
\newcommand{\lb}[1]{\left ( #1 \right )}

\begin{document}
    \title{A Micro-Macro Decomposition-Based Asymptotic-Preserving Random
    Feature Method for Multiscale Radiative Transfer Equations}

    \author[1,2]{Jingrun Chen}
    \author[3,4]{Zheng Ma}
    \author[1, 2, \thanks{Corresponding author: wukekever@ustc.edu.cn}]{Keke Wu}
    
    \affil[1]{School of Mathematical Sciences and Suzhou Institute for Advanced Research, University of Science and Technology of China, Jiangsu,
        215217, P. R. China.}
    \affil[2]{Suzhou Big Data \& AI Research and Engineering Center}
    \affil[3]{School of Mathematical Sciences, Shanghai Jiao Tong University, Shanghai,
    200240, P. R. China.}
    \affil[4]{Institute of Natural Sciences, MOE-LSC,
        Shanghai Jiao Tong University, Shanghai, 200240, P. R. China.}
        
    \date{\today}

    \maketitle

    \begin{abstract}
        This paper introduces the Asymptotic-Preserving Random Feature Method (APRFM) 
        for the efficient resolution of multiscale radiative transfer equations.
        The APRFM effectively addresses the challenges posed by stiffness and
        multiscale characteristics inherent in radiative transfer equations 
        through the application of a micro-macro decomposition strategy. 
        This approach decomposes the distribution
        function into equilibrium and non-equilibrium components, allowing for
        the approximation of both parts through the random feature method (RFM) within a
        least squares minimization framework. 
        The proposed method exhibits remarkable robustness across different scales and achieves high accuracy with fewer degrees of freedom and collocation points than the vanilla RFM. Additionally, compared to the deep neural network-based method, our approach offers significant advantages in terms of parameter efficiency and computational speed. These benefits have been substantiated through numerous numerical experiments conducted on both one- and two-dimensional problems.
    \end{abstract}

    \section{Introduction}

    The radiative transfer equation (RTE) is the governing equation that models
    the propagation and interactions of radiation or particles within
    participating media~\cite{modest2021radiative,larsen1987asymptotic}. It is a fundamental integro-differential equation in various fields, including
    astrophysics, radiative transfer~\cite{marshak20053d}, neutron transport~\cite{lewis1984}
    and optical tomography~\cite{klose2002optical1,klose2002optical2,ren2006frequency,arridge2009optical},
    etc. In recent years, there has been significant interest in devising accurate
    and efficient methods for solving the multiscale radiative transfer equation.
    The primary bottleneck in numerically resolving the radiative transfer
    equation stems from the high dimensionality of phase space, stiffness of collision
    terms, and multiscale features, among others. Various numerical methods have been developed in the field of computational
    methods for the radiative transfer equation, which can generally be categorized into two classes: deterministic
    methods and stochastic simulation methods. One of the most popular
    deterministic methods is the discrete ordinates/velocity method (DOM/DVM)~\cite{lewis1984,adams2002fast},
    sometimes referred to as the $S_{N}$ method. The DOM discretizes the angular
    variable and solves the RTE along the discrete directions. Spherical
    harmonics methods possess the advantage of rotational invariance and are
    widely used in solving the radiative transfer equation~\cite{marshak1947note,evans1998spherical}.
    Note that the scale parameter (Knudsen number) in the radiative transfer
    equation can vary significantly, ranging from the kinetic regime to the diffusive
    regime, therefore numerical methods should be able to handle the multiscale
    nature of the radiative transfer equation. To tackle this challenge, two
    major categories of methods have been developed, namely, the domain
    decomposition-based methods~\cite{golse2003domain} and the asymptotic-preserving
    schemes~\cite{jin2010review}. Domain decomposition-based methods decompose the
    domain into different regions, where different differential equations are
    solved in each region with an interface condition to couple them. The asymptotic-preserving
    schemes aim to design numerical methods that are uniformly stable and accurate,
    regardless of the scale parameter. For stochastic simulation methods, the
    direct simulation Monte Carlo (DSMC) method is widely used for solving
    the radiative transfer equation~\cite{alexander1997direct}. There are also
    some references on solving the radiative transfer equation, interested readers
    can refer to~\cite{frank2007approximate,li2017implicit, shi2023efficient,xiong2022high,liu2023implicit,han2014two,jin2000uniformly}
    for more details.

    In the past few years, deep learning methods have shown great potential in
    solving high-dimensional and complex PDE problems, due to their strong fitting
    ability and generalization ability. Researchers have developed some deep
    learning-based methods for solving the radiative transfer equation and
    kinetic equations~\cite{lou2021physics,CLM,HJJL,lee2021model,XIAO2023112317,
    lu2022solving,wuAPNN,wuAPNNv2,wu2023asymptotic,li2022model,zhang2023asymptotic,li2023learning,
    li2024solving,li2024macroscopic,xu2023transfer,lee2024structure,zhang2024ap}.
    Although deep learning methods have achieved great success in solving PDEs,
    they still face challenges in terms of interpretability, generalization, and
    computational efficiency. For example, the training time for deep learning
    methods can be very long, and the accuracy of the solution may be less than satisfactory.
    The random feature method (RFM) bridges the gap between deep learning methods and
    traditional numerical methods, and is effective in solving
    various PDEs~\cite{chen2022bridging,chi2024random,chen2023random,chen2024optimization}.
    The RFM approximates the solution of the PDE by a linear
    combination of random feature functions and Partition of Unity (PoU) functions,
    and the coefficients are determined by a least squares minimization problem.
    In particular, the random feature functions are constructed by two-layer
    neural networks with fixed parameters in the hidden layers and the PoU
    functions are constructed by tensor product of univariate functions. Benefiting
    from the construction of the random feature functions and least squares minimization
    problem, the RFM can achieve high accuracy and efficiency
    in solving PDEs.

    In this paper, we start by introducing the RFM for solving
    the radiative transfer equation and find that the vanilla random feature
    method has tremendous difficulty in resolving small scales. To address this issue,
    we propose the Asymptotic-Preserving Random Feature Method (APRFM) for solving
    the multiscale radiative transfer equation. The APRFM is designed to
    effectively handle the multiscale nature of the radiative transfer equation
    by utilizing a micro-macro decomposition approach. Our method can approximate
    the solution of the radiative transfer equation by decomposing the distribution
    function into equilibrium and non-equilibrium components and approximating
    both parts through the RFM within a least squares minimization
    framework. The proposed method demonstrates superior robustness across varying
    scales compared to the vanilla RFM and is more efficient
    than the previous deep learning methods.

    The rest of the paper is organized as follows. In Section 2, we gave a
    brief introduction to the radiative transfer equation and the random
    feature method. Besides, we demonstrate the difficulty of the vanilla random
    feature method in resolving small scales. Section 3 is the main
    part of the paper, where we propose the Asymptotic-Preserving Random Feature
    Method. In Section 4, we present numerical results for both one- and two-dimensional
    problems to validate the effectiveness of our method. Finally, we conclude
    the paper in the last section.

    \section{Preliminaries}
    
    \subsection{The radiative transfer equation}

    We consider the scaled form of the stationary radiative transfer equation on
    a bounded Lipschitz domain $D \subset \R^{d}$ as follows:
     \begin{equation}\label{eqn:rte}
        \bv \cdot \nabla_{\bx}f(\bx, \bv) = \frac{\sigma_{s}(\bx)}{\eps}
        \LO f(\bx, \bv) - \eps \sigma_{a}(\bx) f(\bx, \bv) + \eps Q(\bx, \bv), \;
        (\bx,\bv) \in D \times \sS^{d-1},
    \end{equation}
    where $f(\bx,\bv)$ denotes the distribution function of particles at space
    position $\bx \in D$ and velocity direction $\bv \in \sS^{d-1}$, and $Q(\bx,
    \bv)$ represents the source function. The non-negative functions
    $\sigma_{s}(\bx)$ and $\sigma_{a}(\bx)$ correspond to the scattering
    coefficient and absorption coefficient, respectively. The dimensionless
    parameter $\eps > 0$, referred to as the Knudsen number, captures the
    ratio of the mean free path to the characteristic length of the domain. In the
    context of the multiscale problem, the Knudsen number $\eps$ spans
    magnitudes from $\mathcal{O}(1)$ (kinetic regime) to $\eps \ll 1$ (diffusive
    regime). A smaller Knudsen number $\eps$ implies more frequent collisions between
    particles. The operator $\LO$ is defined as:
    \begin{equation}\label{eqn:operator}
        \LO f = \frac{1}{|\sS^{d-1}|}\int_{\sS^{d-1}}k(\bv, \bv
        ') ( f(\bv')-f(\bv)) \diff \bv',
    \end{equation}
    where $k:\sS^{d-1}\times \sS^{d-1}\to \R$ is a non-negative kernel. Here, for
    convenience, we employ the notation
    \begin{equation}\label{eqn:average}
        \avg{h}:= \frac{1}{|\sS^{d-1}|}\int_{\sS^{d-1}}h(\bv')
        \diff \bv',
    \end{equation}
    which denotes the velocity angular average of the function $h$ over the unit
    sphere $\sS^{d-1}$. Next, we make some important assumptions~\cite{egger2014lp}
    on the operator $\LO$ for the well-posedness of the radiative transfer equation~\eqref{eqn:rte}:
    \begin{itemize}
        \item $\avg{\LO f}= 0, \; \forall f(\bv) \in L^{2}(\sS^{d-1})$.

        \item $\LO$ is a self-adjoint and non-positive operator in $L^{2}(\sS^{d-1}
            )$.

        \item The null space of $\LO$ is $\{f \in L^{2}(D \times \sS^{d-1}): f =
            \avg{f}\}$.
    \end{itemize}

    In this paper, we focus on the inflow boundary:
    \begin{equation}\label{eqn:bdy}
        \Gamma_{-}:= \{(\bx,\bv) \in \partial D \times \sS^{d-1}|
        \ \bv \cdot \bn(\bx) < 0 \},
    \end{equation}
    where $\bn(\bx)$ represents the outward normal vector to the boundary
    $\partial D$.

    In the diffusive limit, the radiative transfer equation~\eqref{eqn:rte} can
    be approximated by the elliptic equation~\cite{bardos1984diffusion,lu2022solving}:
    \begin{equation}\label{eqn:elliptic}
        \avg{\bv \cdot \nabla_{\bx} \LO^{-1} \lb{\frac{1}{\sigma_{s}(\bx)} \bv \cdot \nabla_{\bx} \rho(\bx)}}
        = - \sigma_{a}(\bx) \rho(\bx) + \avg{Q(\bx,\bv)},
    \end{equation}
    where $\rho(\bx) = \avg{f(\bx,\bv)}$ denotes the macroscopic density
    function.

    \subsection{Random Feature Method}

    Before proposing the method in this paper, we will give a brief introduction
    to the RFM. Consider the following PDE in a bounded domain
    $\Omega \subset \R^{d}$:
    \begin{equation}\label{eqn:general-pde}
        \mathcal{A}u(\by) = 0, \; \by \in \Omega,
    \end{equation}
    where $\mathcal{A}$ is a differential operator.

    In general, the RFM represents the solution $u(\by)$ by a two-layer
    neural networks with the inner parameters held fixed and chosen randomly.
    First, we introduce a hypercube $\Omega_{c}= \prod_{i = 1}^{d}[a_{i}, b_{i}]$
    of proper size to completely enclose the domain $\Omega$. Here, the notation
    $\prod_{i = 1}^{d}[a_{i}, b_{i}]$ denotes the Cartesian product of $d$ intervals
    $[a_{i}, b_{i}]$, that is, $[a_{1}, b_{1}] \times [a_{2}, b_{2}] \times \cdots
    \times [a_{d}, b_{d}]$. Then we partition the hypercube $\Omega_{c}$ into
    $M$ non-overlapping hyper-rectangles $\Omega_{i}$:
    \begin{equation}\label{eqn:partition}
        \Omega_{c}= \bigcup_{i=1}^{M}\Omega_{i}, \; \Omega_{i} = \prod_{j = 1}^{d}[a_{ij}, b_{ij}],
    \end{equation}
    {and $a_{{i+1},j} = b_{ij}$,}
    which is analogous to the finite element method and one can generalize the partition
    to a more general shape if necessary. Assume
    $\by = (y_{1}, \cdots, y_{d})^{T}\in \Omega$, we apply the following transformation
    to obtain a partition dependent normalized vector $\tilde{\by}_{i}= (\tilde{y}_{i1}, \cdots, \tilde{y}_{id})^{T}$:
    \begin{equation}\label{eqn:transform}
        \tilde{y}_{ij}= 2 \frac{y_{j}- a_{ij}}{b_{ij}- a_{ij}} - 1, \; i = 1, \cdots, M, \; j = 1, \cdots, d.
    \end{equation}
    Denote the center and radius of the hyper-rectangle $\Omega_{i}$ by:
    \begin{equation}\label{eqn:center-raduis}
        \bmu_{i}= \lb{\frac{b_{i1}+ a_{i1}}{2}, \cdots, \frac{b_{id}+ a_{id}}{2}}^{T}, \; 
        \bsigma_{i}= \lb{\frac{b_{i1}- a_{i1}}{2}, \cdots, \frac{b_{id}- a_{id}}{2}}^{T},
    \end{equation}
    then the transformation can be rewritten in the following vector form:
    \begin{equation}\label{eqn:vec-transform}
        \tilde{\by}_{i}= \frac{\by - \bmu_{i}}{\bsigma_{i}}, \; i = 1, \cdots, M.
    \end{equation}
    We can find that the above transformation maps the subdomain $\Omega_{i}$ to
    the hypercube $[-1, 1]^{d}$.

    The RFM consists of two parts: the PoU functions and the random
    feature functions. The PoU functions $\{\psi_{i}(\by)\}_{i=1}^{M}$ is
    construct through the tensor product of univariate function:
    \begin{equation}\label{eqn:pou-fn}
        \psi_{i}(\by) = \prod_{j=1}^{d}\varphi(\tilde{y}_{ij}), \; i = 1, \cdots, M,
    \end{equation}
    where $\varphi$ is chosen in~\cite{chen2022bridging} as:
    \begin{equation}\label{eqn:psi-fn}
        \varphi_{a}(z) =
        \begin{cases}
            1, & |z| \le 1,   \\
            0, & \text{else},
        \end{cases}
        \; \text{or}\; \varphi_{b}(z) =
        \begin{cases}
            1,                                          & |z| \le \frac{3}{4},                \\
            \displaystyle \frac{1 - \sin(2\pi |z|)}{2}, & \frac{3}{4}\le |z| \le \frac{5}{4}, \\
            0,                                          & \text{else}.
        \end{cases}
    \end{equation}
    The random feature functions
    $\{\phi_{ij}(\by)\}_{1 \le i \le M, \; 1 \le j \le J_i}$ are constructed by
    the following formula:
    \begin{equation}\label{eqn:random-feature}
        \phi_{ij}(\by) = \sigma \lb{\bw_{ij} \cdot \tilde{\by}_i + b_{ij}}, \; i = 1, \cdots, M, \; j = 1, \cdots, J_{i},
    \end{equation}
    where $\bw_{ij}\in \R^{d}$ and $b_{ij}\in \R$ are randomly generated with the
    uniform distribution in $[-B, B]$ and fixed parameters, and $\sigma$ is a
    scalar function, named activation function.
    {For simplicity, we assume that all $J_i$ are identical and denote them collectively as $J$.
    Generally, random feature functions are globally defined, whereas PDE solutions often exhibit local variations. 
    To address this, RFM constructs multiple local solutions, each associated with a random feature model, and seamlessly integrates them using a partition of unity (PoU). The introduction of PoU generates local random feature functions, offering a more flexible and general strategy than domain decomposition or mesh generation. The number of partitions $M$ can be regarded as a mechanism for adapting to the solution's local variations. The parameter $B$  is used to control the initialization range of the weights $\{\bw_{ij}\}$ and $\{b_{ij}\}$.}

    Finally, the RFM approximates the solution $u(\by)$
    through the linear combination of random feature functions together with PoU
    functions:
    \begin{equation}\label{eqn:rfm-fn}
        u(\by) \approx u_{M}(\by) = \sum_{i=1}^{M}\psi_{i}(\by)
        \sum_{j=1}^{{J}}u_{ij}\phi_{ij}(\by),
    \end{equation}
    where $\{u_{ij}\}_{1 \le i \le M, \; 1 \le j \le {J}}$ are the unknown coefficients
    to be determined. For simplicity, we denote the set of all coefficients 
    $\{u_{ij}\}_{1 \le i \le M, \; 1 \le j \le {J}}$ by $\theta$. The degrees of freedom
    of the RFM is ${Z = M J}$.

    The problem of determining the coefficients
    $\{u_{ij}\}_{1 \le i \le M, \; 1 \le j \le {J}}$ can be formulated as a
    least squares problem:
    \begin{equation}\label{eqn:least-square}
        \min_{\theta}\| \mathcal{A}u_{M}(\by) \|_{2}^{2},
    \end{equation}
    where $\| \cdot \|_{2}$ denotes the $\ell^{2}$ norm. Besides, the boundary and/or
    initial conditions which are denoted by $\mathcal{B}u_{M}(\by) = 0$ can be incorporated
    into the least squares problem~\eqref{eqn:least-square} as constraints.

    To solve the least squares problem~\eqref{eqn:least-square}, it is necessary
    to generate several collocation points for both $\Omega$ and
    $\partial \Omega$, respectively. Assume the interior collocation points are
    $\{\by_{\text{int}}^{k}\}_{1 \le k \le N_{\text{int}}}$, and the boundary and/or
    initial collocation points are $\{\by_{\text{bdy}}^{k}\}_{1 \le k \le N_{\text{bdy}}}$.
    Then the discrete least squares problem~\eqref{eqn:least-square} can be written
    as:
    \begin{equation}\label{eqn:least-square-discrete}
        \min_{\theta}\sum_{k=1}^{N_{\text{int}}} \lambda_{\text{int}}^k | \mathcal{A} u_{M}(\by_{\text{int}}^{k}) |^{2} 
        + \sum_{k=1}^{N_{\text{bdy}}} \lambda_{\text{bdy}}^k | \mathcal{B}u_{M}(\by_{\text{bdy}}^{k}) |^{2}
        ,
    \end{equation}
    where $\{\lambda_{\text{int}}^k\}$ and $\{\lambda_{\text{bdy}}^k\}$ are regularization
    parameters. The setting of the regularization parameters can be found in~\cite{chen2022bridging}.
    When the operators $\mathcal{A}$ and $\mathcal{B}$ are linear, the minimization
    problem~\eqref{eqn:least-square-discrete} can be solved by the standard
    linear least squares method.

\begin{figure}[htbp]
        \centering
        \includegraphics[width=0.8\textwidth]{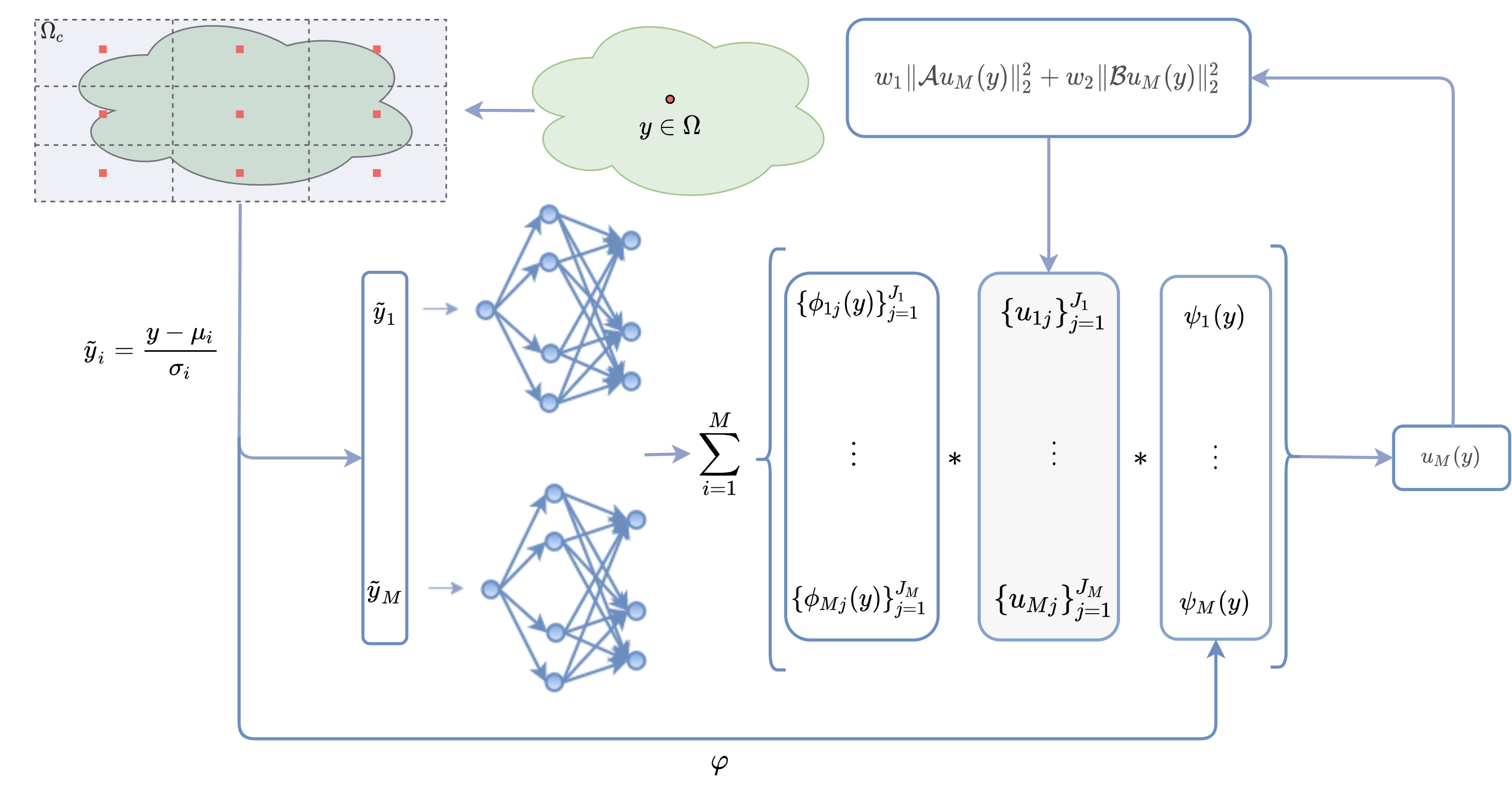}
        \caption{Schematic diagram of the random feature method.}
        \label{fig:rfm}
    \end{figure}

    We present a schematic diagram of the RFM in Figure~\ref{fig:rfm}.

    \subsection{Difficulty of RFM to resolve small scales}\label{sec:motivation}

    In this section, we demonstrate that the vanilla RFM has
    tremendous difficulty in resolving small scales in the multiscale radiative
    transfer equation through a simple example. Consider the following one-dimensional
    radiative transfer equation studied in~\cite{lu2022solving}:
    \begin{equation}\label{eqn:rte1d-motivation}
        \begin{cases}
            \displaystyle \eps v \cdot \nabla_{x}f = \avg{f}- f - \eps v, \; (x, v) \in [0, 1] \times [-1, 1], \\
            f(0, v > 0) = 1, \; f(1, v < 0) = 0.
        \end{cases}
    \end{equation}
    The exact solution of the above equation is $f_{\text{ex}}(x, v) = 1 - x$. 
    We apply the RFM to solve the radiative transfer equation~\eqref{eqn:rte1d-motivation},
    that is, we approximate the solution $f(x, v)$ by the following formula:
    \begin{equation}\label{eqn:rfm-f}
        f(x, v) \approx f_{M}(x, v) = \sum_{i=1}^{M}\psi_{i}(x,v) 
        \sum_{j=1}^{{J}}f_{ij}\phi_{ij}(x, v).
    \end{equation}
    The discrete least squares problem for determining the coefficients $\{f_{ij}\}_{1 \le i \le M, \; 1 \le j \le J_i}$ can be formulated as:
    \begin{equation}\label{eqn:problem-discrete}
        \min_{\{f_{ij}\}_{1 \le i \le M, \; 1 \le j \le J_i}} \mathcal{R}^{\eps} 
        := \mathcal{R}_{\text{int}}^{\eps} + \mathcal{R}_{\text{bdy}}^{\eps},
    \end{equation}
    where $\mathcal{R}_{\text{int}}^{\eps}$ and
    $\mathcal{R}_{\text{bdy}}^{\eps}$ are defined as:
    \begin{equation}\label{eqn:residual-int}
        \mathcal{R}_{\text{int}}^{\eps} = \sum_{k=1}^{N_{\text{int}}} \lambda_{\text{int}}^k
        |\eps v_{\text{int}}^{k} \cdot \nabla_{x} f_{M}(x_{\text{int}}^{k}, v_{\text{int}}^{k}) 
        - \avg{f_{M}}(x_{\text{int}}^{k}) + f_{M}(x_{\text{int}}^{k}, v_{\text{int}}^{k}) 
        + \eps v_{\text{int}}^{k}|^{2},
    \end{equation}
    and
    \begin{equation}\label{eqn:residual-bdy}
        \mathcal{R}_{\text{bdy}}^{\eps} = \sum_{k=1}^{N_{\text{bdy}}} \lambda_{\text{bdy}}^k
        | f_{M}(x_{\text{bdy}}^{k}, v_{\text{bdy}}^{k}) - f_{\text{bdy}}^{k}|^{2},
    \end{equation}
    respectively. 
    Here, $\{(x_{\text{int}}^{k}, v_{\text{int}}^{k})\}_{1 \le k \le N_{\text{int}}}$ 
    and $\{(x_{\text{bdy}}^{k}, v_{\text{bdy}}^{k})\}_{1 \le k \le N_{\text{bdy}}}$ are 
    the interior and boundary collocation points,
    and $\{f_{\text{bdy}}^{k}\}_{1 \le k \le N_{\text{bdy}}}$ are the boundary
    values. The regularization parameters $\{\lambda_{\text{int}}^k\}$ and $\{\lambda_{\text{bdy}}^k\}$
    are positive constants.

    Denote the linear operator
    \begin{equation}\label{eqn:operator-T}
        \mathcal{T}^{\eps}f(x, v) = \eps v \cdot \nabla_{x} f(x, v) - \avg{f}(x) + f(x, v),
    \end{equation}
    we have
    \begin{equation}\label{eqn:T-f_M}
        \mathcal{T}^{\eps}f_{M} = \mathcal{T}^{\eps}\lb{\sum_{i=1}^{M} \psi_i(x, v) \sum_{j=1}^{{J}} f_{ij} \phi_{ij}(x, v)}
        = \sum_{i=1}^{M}\sum_{j=1}^{{J}}f_{ij}\mathcal{T}^{\eps}(\psi_{i}\phi_{ij}).
    \end{equation}

    Introduce the notation
    $s^{\text{int}, k}_{ij} := \mathcal{T}^{\eps}(\psi_{i}\phi_{ij})(x_{\text{int}}
    ^{k}, v_{\text{int}}^{k})$
    and
    $s^{\text{bdy}, k}_{ij} := \psi_{i}(x_{\text{bdy}}^{k}, v_{\text{bdy}}^{k}) 
    \phi_{ij}(x_{\text{bdy}}^{k}, v_{\text{bdy}}^{k})$, we can defined the matrix $\bA$
    and the vector $\bm{b}$
    \begin{equation}
        \bA =
        \begin{bmatrix}
            \bA_{\text{int}} \\
            \bA_{\text{bdy}}
        \end{bmatrix}
        \in \R^{N \times {Z}}, \; \bm{b} =
        \begin{bmatrix}
            \bm{b}_{\text{int}} \\
            \bm{b}_{\text{bdy}}
        \end{bmatrix}
        \in \R^{N}, \; N = N_{\text{int}}+ N_{\text{bdy}},
    \end{equation}
    with
    \begin{equation}\label{eqn:residual-int-matrix}
        \bA_{\text{int}} =
        \begin{bmatrix}
            s^{\text{int}, 1}_{11}              & \cdots & s^{\text{int}, 1}_{1 {J}}              & s^{1}_{21}                          & \cdots & s^{\text{int}, 1}_{M {J}}              \\
            s^{\text{int}, 2}_{11}              & \cdots & s^{\text{int}, 2}_{1 {J}}              & s^{2}_{21}                          & \cdots & s^{\text{int}, 2}_{M {J}}              \\
            \vdots                              & \ddots & \vdots                                 & \vdots                              & \ddots & \vdots                                 \\
            s^{\text{int}, N_{\text{int}}}_{11} & \cdots & s^{\text{int}, N_{\text{int}}}_{1 {J}} & s^{\text{int}, N_{\text{int}}}_{21} & \cdots & s^{\text{int}, N_{\text{int}}}_{M {J}}
        \end{bmatrix}, \, 
        \bm{b}_{\text{int}} = -\eps
        \begin{bmatrix}
            v_{\text{int}}^{1}              \\
            v_{\text{int}}^{2}              \\
            \vdots                          \\
            v_{\text{int}}^{N_{\text{int}}}
        \end{bmatrix},
    \end{equation}
    \begin{equation}\label{eqn:residual-bdy-matrix}
        \bA_{\text{bdy}} =
        \begin{bmatrix}
            s^{\text{bdy}, 1}_{11}              & \cdots & s^{\text{bdy}, 1}_{1 {J}}              & s^{\text{bdy}, 1}_{21}              & \cdots & s^{\text{bdy}, 1}_{M {J}}              \\
            s^{\text{bdy}, 2}_{11}              & \cdots & s^{\text{bdy}, 2}_{1 {J}}              & s^{\text{bdy}, 2}_{21}              & \cdots & s^{\text{bdy}, 2}_{M {J}}              \\
            \vdots                              & \ddots & \vdots                                 & \vdots                              & \ddots & \vdots                                 \\
            s^{\text{bdy}, N_{\text{bdy}}}_{11} & \cdots & s^{\text{bdy}, N_{\text{bdy}}}_{1 {J}} & s^{\text{bdy}, N_{\text{bdy}}}_{21} & \cdots & s^{\text{bdy}, N_{\text{bdy}}}_{M {J}}
        \end{bmatrix}, \; 
        \bm{b}_{\text{bdy}} =
        \begin{bmatrix}
            f_{\text{bdy}}^{1}              \\
            f_{\text{bdy}}^{2}              \\
            \vdots                          \\
            f_{\text{bdy}}^{N_{\text{bdy}}}
        \end{bmatrix}.
    \end{equation}
    Let the vector of coefficients be
    \begin{equation}\label{eqn:rfm-f-vector}
        \bm{f}=
        \begin{bmatrix}
            f_{11} & \cdots & f_{1 {J}} & f_{21} & \cdots & f_{M {J}}
        \end{bmatrix}^{T},
    \end{equation}
    the optimal coefficients $\bm{f}$ of the problem \eqref{eqn:problem-discrete}
    can be obtained by
    \begin{equation}\label{eqn:rfm-problem}
        \bm{f}^{*}= \min_{\bm{f}}\| \bA \bm{f}- \bm{b}\|_{2}^{2}.
    \end{equation}
    The solution $\bm{f}^{*}$ can be obtained by the standard linear least
    squares method.

    First, we consider the problem~\eqref{eqn:rte1d-motivation} with $\eps = 1$
    to validate the effectiveness of the RFM. Here, we set the
    number of partitions $M = 1$ with the PoU function {$\varphi_{b}$}. We choose the
    $\tanh$ function as our activation function and the number of random feature functions
    {$J = 2^{n}, n = 3, 4, 5, 6, 7$}. 
    The weights $\bw_{ij}$ and biases $b_{ij}$ are randomly generated with $B = 1$. 
    Besides, the collocation points are uniform grids with $(N_{x}, N_{v}) = (32, 64)$. The integration of the operator
    $\LO$ is approximated by Gauss-Legendre quadrature with $16$ points.

\begin{figure}[htbp]
        \centering
        \includegraphics[width=0.5\textwidth]{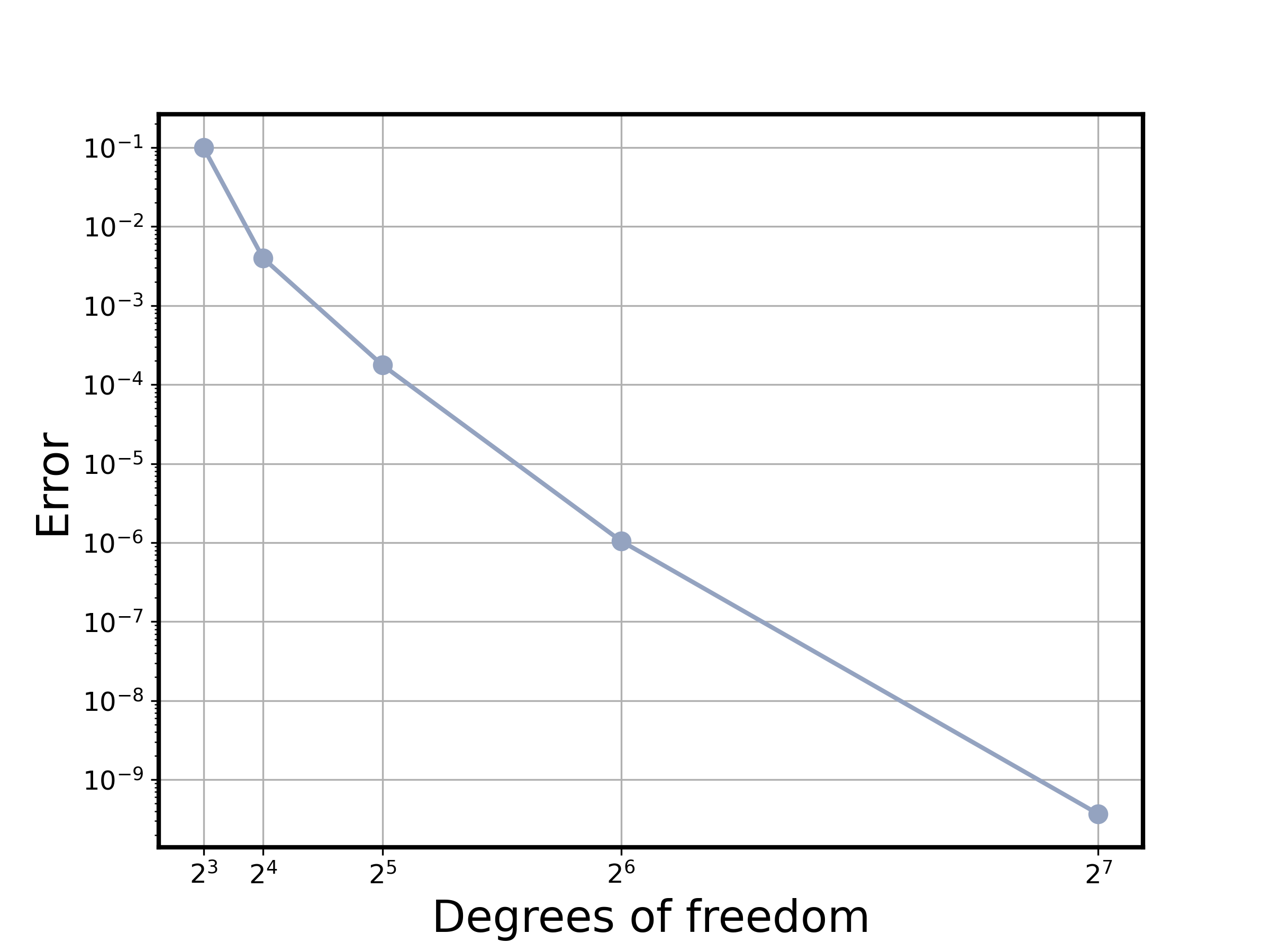}
        \caption{Relative $\ell^{2}$ error of the RFM solution $f_{M}$ with
        respect to the degrees of freedom.}
        \label{fig:rfm-error-semilog}
    \end{figure}

\begin{figure}[htbp]
        \centering
        \includegraphics[width=\textwidth]{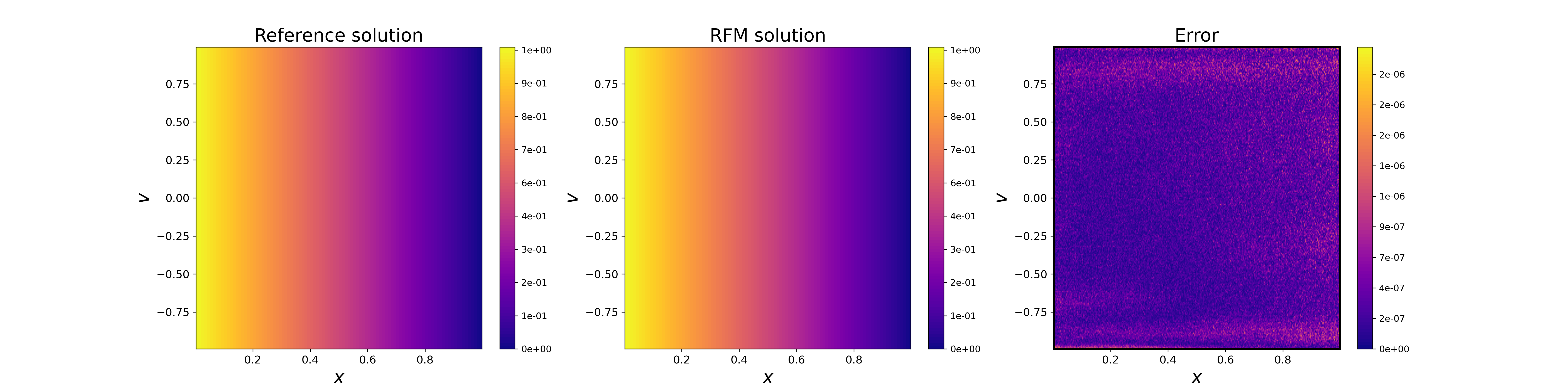}
        \caption{Reference solution v.s. RFM solution.}
        \label{fig:rte1d-1e0}
    \end{figure}

    The relative $\ell^{2}$ error of the solution $f_{M}$ {(as shown in \eqref{eqn:error})} concerning the degrees of freedom {$Z$} is shown in Figure~\ref{fig:rfm-error-semilog}. One can
    observe that the error decreases exponentially to the degrees
    of freedom {$Z$}. We also plot the solution $f_{M}$ obtained by the random feature
    method with {$J = 2^{7}$} and reference solution $f_{\text{ex}}$ in Figure~\ref{fig:rte1d-1e0}.
    The relative $\ell^{2}$ error of the RFM is
    {$3.69 \times 10^{-10}$}.

    Next, we test a sequential of problems~\eqref{eqn:rte1d-motivation} with
    $\eps = 10^{-2^i}$, $i = 1, 2, 3, 4$ with the same setting as the previous problem
    except for the number of random features ${J}$ and the collocation points $(N_{x}
    , N_{v})$. We record the relative $\ell^{2}$ error of the RFM solution to the degrees of freedom and the number of collocation points in
    Table~\ref{tab:dofs-dependency-rfm} and Table~\ref{tab:sample-dependency-rfm}, respectively. 
    One can find that the small-scale problems ($\eps \ll 1$) are much more difficult to resolve than the
    large-scale problems in terms of the degrees of freedom and the number of collocation points.
    
    \begin{table}[tbhp]
        \caption{The dependency of the RFM on degrees of freedom. The collocation
        points are uniform grids with $(N_{x}, N_{v}) = (64, 128)$.}
        \label{tab:dofs-dependency-rfm}
        \centering
        \begin{tabular}{ccccccc}
            \toprule[1pt] \noalign{\smallskip} \multirow{2}*{$\eps$}                           & \multicolumn{5}{c}{$J$} \\
                                                                                               & \multicolumn{1}{c}{$16$}   & \multicolumn{1}{c}{$32$} & \multicolumn{1}{c}{$64$} & \multicolumn{1}{c}{$128$} & \multicolumn{1}{c}{$256$} \\
            \noalign{\smallskip} \midrule[1pt] \noalign{\smallskip} \multirow{1}*{{$10^{-2}$}} 
            & ${\text{1.87 e-}1}$        & ${\text{8.92 e-}3}$      & ${\text{4.95 e-}5}$      & ${\text{5.59 e-}9}$       & ${\text{5.76 e-}11}$       \\
            \multirow{1}*{{$10^{-4}$}}   
            & ${\text{2.10 e-}1}$        & ${\text{6.88 e-}2}$      & ${\text{1.17 e-}2}$      & ${\text{1.93 e-}5}$       & ${\text{2.92 e-}7}$       \\
            \multirow{1}*{{$10^{-8}$}}   
            & ${\text{1.83 e-}1}$        & ${\text{6.92 e-}2}$      & ${\text{6.08 e-}2}$      & ${\text{1.67 e-}3}$       & ${\text{2.54 e-}3}$       \\
            \multirow{1}*{{$10^{-16}$}}  
            & ${\text{1.54 e-}1}$        & ${\text{1.39 e-}1}$      & ${\text{6.06 e-}2}$      & ${\text{3.01 e-}2}$       & ${\text{6.02 e-}2}$       \\
            \noalign{\smallskip} \bottomrule[1pt]
        \end{tabular}
    \end{table}
    
    \begin{table}[tbhp]
        \caption{The dependency of the RFM on the number of collocation points. The
        number of random feature functions $J= 128$.}
        \label{tab:sample-dependency-rfm}
        \centering
        \begin{tabular}{cccccc}
            \toprule[1pt] \noalign{\smallskip} \multirow{2}*{$\eps$}                           & \multicolumn{4}{c}{$(N_{x}, \, N_{v})$} \\
            & \multicolumn{1}{c}{$(16, \, 32)$}      & \multicolumn{1}{c}{$(32, \, 64)$} & \multicolumn{1}{c}{$(64, \, 128)$} & \multicolumn{1}{c}{$(128, \, 256)$} \\
            \noalign{\smallskip} \midrule[1pt] \noalign{\smallskip} \multirow{1}*{{$10^{-2}$}} 
            & ${\text{5.03 e-}8}$                    & ${\text{5.59 e-}9}$               & ${\text{6.94 e-}9}$                & ${\text{3.80 e-}10}$                 \\
            \multirow{1}*{{$10^{-4}$}}                                                         
            & ${\text{1.09 e-}5}$                    & ${\text{1.93 e-}5}$               & ${\text{2.62 e-}5}$                & ${\text{8.98 e-}6}$                 \\
            \multirow{1}*{{$10^{-8}$}}                                                         
            & ${\text{3.64 e-}2}$                    & ${\text{1.67 e-}3}$               & ${\text{1.10 e-}3}$                & ${\text{4.96 e-}4}$                 \\
            \multirow{1}*{{$10^{-16}$}}                                                        
            & ${\text{3.08 e-}2}$                    & ${\text{3.01 e-}2}$               & ${\text{4.01 e-}2}$                & ${\text{3.89 e-}2}$                 \\
            \noalign{\smallskip} \bottomrule[1pt]
        \end{tabular}
    \end{table}
    
    Consider the problem~\eqref{eqn:rte1d-motivation} with small scale, i.e.,
    $\eps \ll 1$. It can be observed that any function $f$ independent of $v$, 
    which satisfies the inflow boundary condition, such as $f = (1 - x)^n$ for $n \ge 2$, 
    will lead to the following:
    \begin{equation}\label{eqn:error-eps}
        \| \eps v \cdot \nabla_{x} f - \avg{f} + f + \eps v \|_2^{2} = \mathcal{O}(\eps^{2}),
    \end{equation}
    but
    \begin{equation}\label{eqn:order-one}
        \| f - f_{\text{ex}}\|_{2}^{2} = \mathcal{O}(1).
    \end{equation}
    {To resolve the multi-scale information, a natural consideration is to increase the number of degrees of freedom and the number of collocation points. The above strategy will lead to a large matrix $\bA$ with a high condition number and make the least squares problem~\eqref{eqn:rfm-problem} ill-conditioned. Moreover, as shown in Table~\ref{tab:sample-dependency-rfm}, when $\eps$ is relatively large, for instance, $\eps = 10^{-2}$, increasing the number of collocation points can substantially enhance the accuracy. However, when $\eps$ is exceedingly small, such as $\eps = 10^{-16}$, increasing the number of collocation points has a negligible impact and fails to significantly improve accuracy. 
    Beyond the aspects previously discussed, we also investigated the impact of partition of unity (PoU) on this problem.
    The domain partitioning is denoted by $(M_{x}, \, M_{v})$. 
    The relative $\ell^{2}$ error are recorded in Table~\ref{tab:pou-dependency-rfm}.
    It can be observed that when $\eps \ll 1$, the effectiveness of the PoU is also quite limited.
    In some cases, it may even result in a decrease in accuracy.
    As can be seen from equation~\eqref{eqn:error-eps}, relying solely on the original radiative transport equation imposes extremely stringent requirements on the mesh size of the partition.}
    
    \begin{table}[tbhp]
        \caption{The dependency of the RFM on PoU. The collocation
        points are uniform grids with $(N_{x}, N_{v}) = (64, 128)$ and the number of random feature functions $J=128$.}
        \label{tab:pou-dependency-rfm}
        \centering
        \begin{tabular}{ccccccc}
            \toprule[1pt] \noalign{\smallskip} \multirow{2}*{$\eps$}                           & \multicolumn{5}{c}{$(M_{x}, \, M_{v})$} \\
                                                                                               & \multicolumn{1}{c}{$(1, \, 1)$}   & \multicolumn{1}{c}{$(2, \, 1)$} & \multicolumn{1}{c}{$(1, \, 2)$} & \multicolumn{1}{c}{$(4, \, 1)$} & \multicolumn{1}{c}{$(1, \, 4)$} \\
            \noalign{\smallskip} \midrule[1pt] \noalign{\smallskip} \multirow{1}*{{$10^{-2}$}} 
            & ${\text{6.94 e-}9}$        & ${\text{1.68 e-}9}$      & ${\text{2.51 e-}9}$      & ${\text{4.26 e-}8}$       & ${\text{8.94 e-}9}$       \\
            \multirow{1}*{{$10^{-4}$}}   
            & ${\text{2.62 e-}5}$        & ${\text{1.55 e-}6}$      & ${\text{5.88 e-}5}$      & ${\text{1.26 e-}5}$       & ${\text{3.18 e-}4}$       \\
            \multirow{1}*{{$10^{-8}$}}   
            & ${\text{1.10 e-}3}$        & ${\text{1.57 e-}2}$      & ${\text{7.83 e-}3}$      & ${\text{1.12 e-}1}$       & ${\text{3.01 e-}2}$       \\
            \multirow{1}*{{$10^{-16}$}}  
            & ${\text{4.01 e-}2}$        & ${\text{6.49 e-}1}$      & ${\text{9.38 e-}2}$      & ${\text{6.07 e-}1}$       & ${\text{8.10 e-}2}$       \\
            \noalign{\smallskip} \bottomrule[1pt]
        \end{tabular}
    \end{table}
    
    The above analysis implies that the
    vanilla RFM based on least squares formulation for the radiative
    transfer equation~\eqref{eqn:rte} has tremendous difficulty in resolving
    small scales in the radiative transfer equation. This is not only due to the
    ill-conditioning of the least squares problem, but also the lack of
    capturing the small-scale part of the solution $f(x, v)$. 
    Such a similar phenomenon
    is also observed in the numerical experiments in~\cite{wuAPNN,lu2022solving}
    when solving the multiscale kinetic equations with small parameters through deep
    neural networks.

    \section{Asymptotic-Preserving Random Feature Method}

\begin{figure}[htbp]
        \centering
        \includegraphics[width=0.5\textwidth]{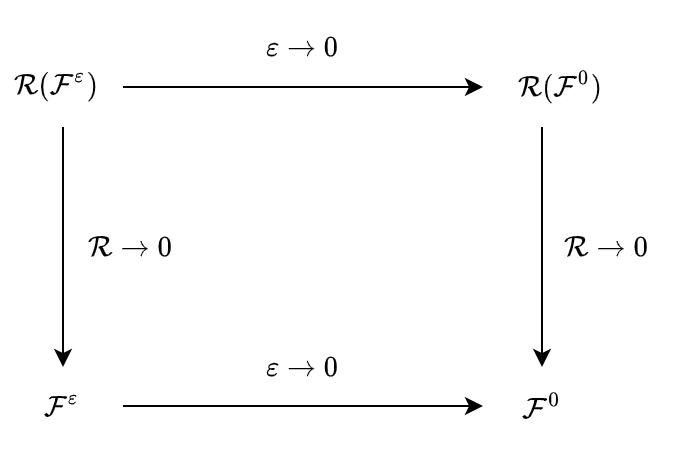}
        \caption{Schematic diagram of the asymptotic-preserving random feature
        method.}
        \label{fig:aprfm}
    \end{figure}

    To address the issue discussed in the previous section, we propose a new random
    feature method based on a micro-macro decomposition~\cite{liu04,lemou2008new},
    called the asymptotic-preserving random feature method. We introduce the
    definition of the asymptotic-preserving random feature method as follows (see
    Figure~\ref{fig:aprfm}):
    {
    \begin{defi}
        Assume $\mathcal{F^{\eps}}$ is the multi-scale model that depends on the
        scale parameter $\eps$ and $\mathcal{F}^{0}$ is the corresponding asymptotic limit model as $\eps \to 0$. Define $\mathcal{R}(\mathcal{F^{\eps}})$ as the least-squares formulation of the model $\mathcal{F^{\eps}}$ when solved using the random feature method. If $\mathcal{R}(\mathcal{F^{\eps}})$ converges to 
        $\mathcal{R}(\mathcal{F}^{0})$ as $\eps \to 0$, and this limit is precisely the least-squares formulation of the limit model $\mathcal{F}^{0}$, then
        the method is called asymptotic-preserving.
    \end{defi}
    }

    Let us revisit the radiative transfer equation~\eqref{eqn:rte},
    \begin{equation}\label{eqn:rte-previous}
        \bv \cdot \nabla_{\bx}f(\bx, \bv) = \frac{\sigma_{s}(\bx)}{\eps}
        \LO f(\bx, \bv) -\eps \sigma_{a}(\bx) f(\bx, \bv) + \eps Q(\bx, \bv),
    \end{equation}
    and we decompose the distribution function $f(\bx, \bv)$ into two parts:
    \begin{equation}\label{eqn:decomposition}
        f(\bx, \bv) = \rho(\bx) + \eps g(\bx, \bv),
    \end{equation}
    where $\rho(\bx) = \avg{f(\bx, \bv)}$ denotes the equilibrium part, and 
    $g(\bx, \bv)$ denotes the non-equilibrium part satisfying $\avg{g(\bx, \bv)}= 0$.

    Substitute the decomposition~\eqref{eqn:decomposition} into the radiative
    transfer equation~\eqref{eqn:rte-previous}, and according to the properties
    of the operator $\LO$, we obtain the following equations:
    \begin{equation}\label{eqn:rte-decomposition}
        \bv \cdot \nabla_{\bx}(\rho(\bx) + \eps g(\bx, \bv)) = 
        \sigma_{s}(\bx)\LO g(\bx, \bv) -\eps \sigma_{a}(\bx) (\rho(\bx)
        + \eps g(\bx, \bv)) + \eps Q(\bx, \bv).
    \end{equation}
    Integrating the above equation over the velocity space $\sS^{d-1}$, we have
    \begin{equation}\label{eqn:macro}
        \avg{\bv \cdot \nabla_{\bx} g}= - \sigma_{a}\rho + \avg{Q}.
    \end{equation}
    Define the orthogonal projection operator
    \begin{equation}\label{eqn:Pi}
        \Pi: \Pi(\cdot)(\bm{v}) = \avg{\cdot},
    \end{equation}
    and the identity operator $\text{Id}$, one can apply the operator
    $\text{Id}- \Pi$ to the equation~\eqref{eqn:rte-decomposition} to obtain the following equation:
    \begin{equation}\label{eqn:micro}
        \bv \cdot \nabla_{\bx}\rho + \eps (\text{Id}- \Pi) (\bv
        \cdot \nabla_{\bx}g) = \sigma_{s}\LO g - \eps^{2}\sigma_{a}g + \eps (Q -
        \avg{Q}).
    \end{equation}
    Therefore, {equations~\eqref{eqn:macro}} and \eqref{eqn:micro} constitute the micro-macro
    system of the radiative transfer equation~\eqref{eqn:rte-previous}:
    \begin{equation}
        \label{eqn:micro-macro}
        \begin{cases}
            \avg{\bv \cdot \nabla_{\bx} g}+ \sigma_{a}\rho = \avg{Q},                                                                                  \\
            \bv \cdot \nabla_{\bx}\rho + \eps (\text{Id}- \Pi) (\bv \cdot \nabla_{\bx}g) 
            = \sigma_{s}\LO g - \eps^{2}\sigma_{a}g + \eps (Q - \avg{Q}). 
        \end{cases}
    \end{equation}

    Our idea is to approximate $\rho(\bx)$ and $g(\bx, \bv)$ by the following
    formula:
    \begin{equation}\label{eqn:aprfm-rho}
        \rho(\bx) \approx \rho_{M}(\bx) = \sum_{i=1}^{M^{\rho}}
        \psi_{i}^{\rho}(\bx) \sum_{j=1}^{{J^{\rho}}}\rho_{ij}\phi_{ij}^{\rho}(\bx
        ),
    \end{equation}
    and
    \begin{equation}\label{eqn:aprfm-g}
        g(\bx, \bv) \approx g_{M}(\bx, \bv) = \sum_{i=1}^{M^{g}}
        \psi_{i}^{g}(\bx, \bv) \sum_{j=1}^{{J^{g}}}g_{ij}\phi_{ij}^{g}(\bx, \bv),
    \end{equation}
    respectively. The coefficients
    $\{\rho_{ij}\}_{1 \le i \le M^{\rho}, \; 1 \le j \le {J^{\rho}}}$ and
    $\{g_{ij}\}_{1 \le i \le M^g, \; 1 \le j \le {J^g}}$ still denoted by
    $\bm{\theta}$ can be determined by the least squares problem of the micro-macro
    system~\eqref{eqn:micro-macro} of the radiative transfer equation~\eqref{eqn:rte}:
    \begin{equation}\label{eqn:ap-residual}
        \min_{\bm{\theta}}\mathcal{R}^{\eps}:= 
        \mathcal{R}_{\text{int}}^{\eps}+ \mathcal{R}_{\text{bdy}}
        ^{\eps},
    \end{equation}
    where $\mathcal{R}_{\text{int}}^{\eps}$ and $\mathcal{R}_{\text{bdy}}^{\eps}$
    are defined as:
    \begin{equation}
        \label{eqn:ap-residual-int}
        \begin{aligned}
            \mathcal{R}_{\text{int}}^{\eps} & = \sum_{k=1}^{N_{\text{int}}} \lambda_{\text{int}, 1}^k | \avg{\bv \cdot \nabla_{\bx} g_M}(\bx_{\text{int}}^{k}) + \sigma_{a}(\bx_{\text{int}}^{k}) \rho_{M}(\bx_{\text{int}}^{k}) - \avg{Q}(\bx_{\text{int}}^{k}) |^{2}                         \\
                                            & + \sum_{k=1}^{N_{\text{int}}} \lambda_{\text{int}, 2}^k | \bv_{\text{int}}^{k}\cdot \nabla_{\bx}\rho_{M}(\bx_{\text{int}}^{k}) + \eps (\text{Id}- \Pi) (\bv_{\text{int}}^{k}\cdot \nabla_{\bx}g_{M}(\bx_{\text{int}}^{k}, \bv_{\text{int}}^{k})) \\
                                            & \qquad \qquad - \sigma_{s}(\bx_{\text{int}}^{k}) \LO g_{M}(\bx_{\text{int}}^{k}, \bv_{\text{int}}^{k}) + \eps^{2}\sigma_{a}(\bx_{\text{int}}^{k}) g_{M}(\bx_{\text{int}}^{k}, \bv_{\text{int}}^{k})                   \\
                                            & \qquad \qquad - \eps (Q(\bx_{\text{int}}^{k}, \bv_{\text{int}}^{k}) - \avg{Q}(\bx_{\text{int}}^{k})) |^{2},
        \end{aligned}
    \end{equation}
    and
    \begin{equation}
        \label{eqn:ap-residual-bdy}\mathcal{R}_{\text{bdy}}^{\eps}= \sum_{k=1}^{N_{\text{bdy}}} \lambda_{\text{bdy}}^k
        | \rho_{M}(\bx_{\text{bdy}}^{k}) + \eps g_{M}(\bx_{\text{bdy}}^{k}, \bv_{\text{bdy}}
        ^{k}) - f_{\text{bdy}}^{k}|^{2},
    \end{equation}
    respectively. Here,
    $\{(\bx_{\text{int}}^{k}, \bv_{\text{int}}^{k})\}_{1 \le k \le N_{\text{int}}}$
    and
    $\{(\bx_{\text{bdy}}^{k}, \bv_{\text{bdy}}^{k})\}_{1 \le k \le N_{\text{bdy}}}$
    are the interior and boundary collocation points, and $\{f_{\text{bdy}}^{k}\}
    _{1 \le k \le N_{\text{bdy}}}$ are the boundary values.
    $\{\lambda_{\text{int}, 1}^k\}, \{\lambda_{\text{int}, 2}^k\}$ and $\{\lambda_{\text{bdy}}^k\}$ are regularization parameters.
    {The setting of the regularization parameters will be determined 
    through a simple but effective way in \eqref{eqn:resacle-factor}.} 
    The degrees of freedom of the APRFM are {$Z = M^{\rho} J^{\rho}+ M^{g}J^{g}$.}

    {
    Taking $\eps \to 0$, the least squares problem~\eqref{eqn:ap-residual-int} reduces to
    \begin{equation}\label{eqn:ap-residual-int-limit}
        \begin{aligned}
            \mathcal{R}_{\text{int}}^{0} & = \sum_{k=1}^{N_{\text{int}}} \lambda_{\text{int}, 1}^k | \avg{\bv \cdot \nabla_{\bx} g_M}(\bx_{\text{int}}^{k}) + \sigma_{a}(\bx_{\text{int}}^{k}) \rho_{M}(\bx_{\text{int}}^{k}) - \avg{Q}(\bx_{\text{int}}^{k}) |^{2}                         \\
                                            & + \sum_{k=1}^{N_{\text{int}}} \lambda_{\text{int}, 2}^k | \bv_{\text{int}}^{k}\cdot \nabla_{\bx}\rho_{M}(\bx_{\text{int}}^{k})  - \sigma_{s}(\bx_{\text{int}}^{k}) \LO g_{M}(\bx_{\text{int}}^{k}, \bv_{\text{int}}^{k}) |^{2},
        \end{aligned}
    \end{equation}
    which corresponds to the least squares formulation of the system
    \begin{equation}\label{eqn:ap-eqn-int-limit}
        \begin{cases}
            \avg{\bv \cdot \nabla_{\bx} g}+ \sigma_{a}\rho = \avg{Q},   \\
            \bv \cdot \nabla_{\bx}\rho = \sigma_{s} \LO g. 
        \end{cases}
    \end{equation}
    From the second equation of~\eqref{eqn:ap-eqn-int-limit}, 
    we obtain $g = \LO^{-1}(\frac{1}{\sigma_{s}(\bx)} \bv \cdot \nabla_{\bx}\rho)$.
    Substituting this into the first equation yields the limiting equation~\eqref{eqn:elliptic}.
    Thus, this proposed method is asymptotic-preserving.
    }

    Denote the linear operators
    \begin{equation}\label{eqn:operator-T11}
        \mathcal{T}_{11}^{\eps}\rho = \sigma_{a}(\bx) \rho
        (\bx),
    \end{equation}
    \begin{equation}\label{eqn:operator-T12}
        \mathcal{T}_{12}^{\eps}g = \avg{\bv \cdot \nabla_{\bx} g(\bx, \bv)},
    \end{equation}
    \begin{equation}\label{eqn:operator-T21}
        \mathcal{T}_{21}^{\eps}\rho = \bv \cdot \nabla_{\bx}
        \rho(\bx),
    \end{equation}
    and
    \begin{equation}\label{eqn:operator-T22}
        \mathcal{T}_{22}^{\eps}g = \eps (\text{Id}- \Pi)
        (\bv \cdot \nabla_{\bx}g(\bx, \bv)) - \sigma_{s}(\bx) \LO g(\bx, \bv) + \eps
        ^{2}\sigma_{a}(\bx) g(\bx, \bv).
    \end{equation}
    We have
    \begin{equation}\label{eqn:T1}
        \mathcal{T}_{11}^{\eps}\rho_{M} = 
        \sum_{i=1}^{M^{\rho}}\sum_{j=1}
        ^{{J^{\rho}}}\rho_{ij}\mathcal{T}_{11}^{\eps}(\psi_{i}^{\rho}\phi_{ij}^{\rho}
        ), \; \mathcal{T}_{12}^{\eps}g_{M}= \sum_{i=1}^{M^{g}}\sum_{j=1}^{{J^{g}}}
        g_{ij}\mathcal{T}_{12}^{\eps}(\psi_{i}^{g}\phi_{ij}^{g}),
    \end{equation}
    and
    \begin{equation}\label{eqn:T2}
        \mathcal{T}_{21}^{\eps}\rho_{M} = 
        \sum_{i=1}^{M^{\rho}}\sum_{j=1}
        ^{{J^{\rho}}}\rho_{ij}\mathcal{T}_{21}^{\eps}(\psi_{i}^{\rho}\phi_{ij}^{\rho}
        ), \; \mathcal{T}_{22}^{\eps}g_{M}= \sum_{i=1}^{M^{g}}\sum_{j=1}^{{J^{g}}}
        g_{ij}\mathcal{T}_{22}^{\eps}(\psi_{i}^{g}\phi_{ij}^{g}).
    \end{equation}

    Introduce the notations
    \begin{equation}\label{eqn:s_11}
        ^{11}s^{\text{int}, k}_{ij} := 
        \mathcal{T}_{11}^{\eps}(\psi_{i}^{\rho}\phi_{ij}^{\rho})(\bx_{\text{int}}^{k}, \bv_{\text{int}}^{k}),
    \end{equation}
    \begin{equation}\label{eqn:s_12}
        ^{12}s^{\text{int}, k}_{ij} := 
        \mathcal{T}_{12}^{\eps}(\psi
        _{i}^{g}\phi_{ij}^{g})(\bx_{\text{int}}^{k}, \bv_{\text{int}}^{k}),
    \end{equation}
    \begin{equation}\label{eqn:s_21}
        ^{21}s^{\text{int}, k}_{ij}:= \mathcal{T}_{21}^{\eps}(\psi_{i}^{\rho}
        \phi_{ij}^{\rho})(\bx_{\text{int}}^{k}, \bv_{\text{int}}^{k}),
    \end{equation}
    \begin{equation}\label{eqn:s_22}
        ^{22}s^{\text{int}, k}_{ij}:= \mathcal{T}_{22}^{\eps}(\psi_{i}^{g}
        \phi_{ij}^{g})(\bx_{\text{int}}^{k}, \bv_{\text{int}}^{k}),
    \end{equation}
    and
    \begin{equation}\label{eqn:s_31}
        ^{31}s^{\text{bdy}, k}_{ij}:= (\psi_{i}^{\rho}\phi_{ij}^{\rho})(\bx_{\text{bdy}}^{k}),
    \end{equation}
    \begin{equation}\label{eqn:s_32}
        ^{32}s^{\text{bdy}, k}_{ij}:= (\eps \psi_{i}^{g}\phi_{ij}^{g})(\bx_{\text{bdy}}^{k}, \bv_{\text{bdy}}^{k}),
    \end{equation}
    then we construct the matrix $\bA$ and the vector $\bm{b}$ as follows:
    \begin{equation}
        \bA =
        \begin{bmatrix}
            \bA_{\text{int}} \\
            \bA_{\text{bdy}}
        \end{bmatrix}
        \in \R^{N \times Z}, \; \bm{b}=
        \begin{bmatrix}
            \bm{b}_{\text{int}} \\
            \bm{b}_{\text{bdy}}
        \end{bmatrix}
        \in \R^{N}, \; N = 2 N_{\text{int}}+ N_{\text{bdy}},
    \end{equation}
    \begin{equation}\label{eqn:ap-residual-matrix-A}
        \bA_{\text{int}}=
        \begin{bmatrix}
            ^{11}\bA_{\text{int}}^{1}              & ^{12}\bA_{\text{int}}^{1}              \\
            ^{21}\bA_{\text{int}}^{1}              & ^{22}\bA_{\text{int}}^{1}              \\
            ^{11}\bA_{\text{int}}^{2}              & ^{12}\bA_{\text{int}}^{2}              \\
            ^{21}\bA_{\text{int}}^{2}              & ^{22}\bA_{\text{int}}^{2}              \\
            \vdots                                 & \vdots                                 \\
            ^{11}\bA_{\text{int}}^{N_{\text{int}}} & ^{12}\bA_{\text{int}}^{N_{\text{int}}} \\
            ^{21}\bA_{\text{int}}^{N_{\text{int}}} & ^{22}\bA_{\text{int}}^{N_{\text{int}}}
        \end{bmatrix}, \; \bA_{\text{bdy}}=
        \begin{bmatrix}
            ^{31}\bA_{\text{bdy}}^{1}              & ^{32}\bA_{\text{bdy}}^{1}              \\
            ^{31}\bA_{\text{bdy}}^{2}              & ^{32}\bA_{\text{bdy}}^{2}              \\
            \vdots                                 & \vdots                                 \\
            ^{31}\bA_{\text{bdy}}^{N_{\text{bdy}}} & ^{32}\bA_{\text{bdy}}^{N_{\text{bdy}}}
        \end{bmatrix},
    \end{equation}
    and
    \begin{equation}\label{eqn:ap-residual-vector-b}
        \bm{b}_{\text{int}}=
        \begin{bmatrix}
            \avg{Q}(\bx_{\text{int}}^{1})                                                                                               \\
            \eps (Q(\bx_{\text{int}}^{1}, \bv_{\text{int}}^{1}) - \avg{Q}(\bx_{\text{int}}^{1}))                                        \\
            \avg{Q}(\bx_{\text{int}}^{2})                                                                                               \\
            \eps (Q(\bx_{\text{int}}^{2}, \bv_{\text{int}}^{2}) - \avg{Q}(\bx_{\text{int}}^{2}))                                        \\
            \vdots                                                                                                                      \\
            \avg{Q}(\bx_{\text{int}}^{N_{\text{int}}})                                                                                  \\
            \eps (Q(\bx_{\text{int}}^{N_{\text{int}}}, \bv_{\text{int}}^{N_{\text{int}}}) - \avg{Q}(\bx_{\text{int}}^{N_{\text{int}}}))
        \end{bmatrix}, \; \bm{b}_{\text{bdy}}=
        \begin{bmatrix}
            f_{\text{bdy}}^{1}              \\
            f_{\text{bdy}}^{2}              \\
            \vdots                          \\
            f_{\text{bdy}}^{N_{\text{bdy}}}
        \end{bmatrix},
    \end{equation}
    where the vectors $^{11}\bA_{\text{int}}^{k},^{12}\bA_{\text{int}}^{k},^{21}\bA
    _{\text{int}}^{k},^{22}\bA_{\text{int}}^{k},^{31}\bA_{\text{bdy}}^{k},^{32}\bA
    _{\text{bdy}}^{k}$ are defined as
    \begin{equation}
        ^{11}\bA_{\text{int}}^{k}=
        \begin{bmatrix}
            ^{11}s^{\text{int}, k}_{11} & \cdots & ^{11}s^{\text{int}, k}_{1 {J^{\rho}}} & ^{11}s^{\text{int}, k}_{21} & \cdots & ^{11}s^{\text{int}, k}_{M^{\rho} {J^{\rho}}}
        \end{bmatrix},
    \end{equation}
    \begin{equation}
        ^{12}\bA_{\text{int}}^{k}=
        \begin{bmatrix}
            ^{12}s^{\text{int}, k}_{11} & \cdots & ^{12}s^{\text{int}, k}_{1 {J^{g}}} & ^{12}s^{\text{int}, k}_{21} & \cdots & ^{12}s^{\text{int}, k}_{M^{g} {J^{g}}}
        \end{bmatrix},
    \end{equation}
    \begin{equation}
        ^{21}\bA_{\text{int}}^{k}=
        \begin{bmatrix}
            ^{21}s^{\text{int}, k}_{11} & \cdots & ^{21}s^{\text{int}, k}_{1 {J^{\rho}}} & ^{21}s^{\text{int}, k}_{21} & \cdots & ^{21}s^{\text{int}, k}_{M^{\rho} {J^{\rho}}}
        \end{bmatrix},
    \end{equation}
    \begin{equation}
        ^{22}\bA_{\text{int}}^{k}=
        \begin{bmatrix}
            ^{22}s^{\text{int}, k}_{11} & \cdots & ^{22}s^{\text{int}, k}_{1 {J^{g}}} & ^{22}s^{\text{int}, k}_{21} & \cdots & ^{22}s^{\text{int}, k}_{M^{g} {J^{g}}}
        \end{bmatrix},
    \end{equation}
    \begin{equation}
        ^{31}\bA_{\text{bdy}}^{k}=
        \begin{bmatrix}
            ^{31}s^{\text{bdy}, k}_{11} & \cdots & ^{31}s^{\text{bdy}, k}_{1 {J^{\rho}}} & ^{31}s^{\text{bdy}, k}_{21} & \cdots & ^{31}s^{\text{bdy}, k}_{M^{\rho} {J^{\rho}}}
        \end{bmatrix},
    \end{equation}
    \begin{equation}
        ^{32}\bA_{\text{bdy}}^{k}=
        \begin{bmatrix}
            ^{32}s^{\text{bdy}, k}_{11} & \cdots & ^{32}s^{\text{bdy}, k}_{1 {J^{g}}} & ^{32}s^{\text{bdy}, k}_{21} & \cdots & ^{32}s^{\text{bdy}, k}_{M^{g} {J^{g}}}
        \end{bmatrix}.
    \end{equation}
    Let $\bm{\theta}= [\bm{\rho}, \bm{g}]^{T}$, with
    \begin{equation}
        \label{eqn:aprfm-rho-vector}\bm{\rho}=
        \begin{bmatrix}
            \rho_{11} & \cdots & \rho_{1 {J^{\rho}}} & \rho_{21} & \cdots & \rho_{M^{\rho} {J^{\rho}}}
        \end{bmatrix},
    \end{equation}
    and
    \begin{equation}\label{eqn:aprfm-g-vector}
        \bm{g}=
        \begin{bmatrix}
            g_{11} & \cdots & g_{1 {J^{g}}} & g_{21} & \cdots & g_{M^{g} {J^{g}}}
        \end{bmatrix},
    \end{equation}
    then the least squares problem \eqref{eqn:ap-residual} can be written as
    \begin{equation}\label{eqn:aprfm-problem}
        \bm{\theta}^{*}= \min_{\bm{\theta}}\| \bA \bm{\theta}
        - \bm{b}\|_{2}^{2}.
    \end{equation}
    
    {To reduce the condition number of the least squares problem, we rescale the residuals of each term in \eqref{eqn:aprfm-problem} based on the maximum absolute value, ensuring that their magnitudes are on the same order. Specifically, the matrix $\bA = \left ( a_{kj} \right ) \in \R^{N \times Z}$, and the vector $\bm{b} = [b_1, \cdots, b_N]^T \in \R^{N}$ are redefined and the regularization parameters $\lambda^k \in \{\lambda_{\text{int}, 1}^k, \lambda_{\text{int}, 2}^k, \lambda_{\text{bdy}}^k\}$ can be defined as follows:
    \begin{equation}\label{eqn:resacle-factor}
        \lambda^k = \frac{1}{\max_{1 \le j \le Z} | a_{kj} |}, \; k = 1, \cdots, N.
    \end{equation}
    Correspondingly, the elements of the matrix $\bA$ and the vector $\bm{b}$ in \eqref{eqn:ap-residual} should be transformed, with their respective elements given by
    \begin{equation}\label{eqn:resacle-elements}
        a_{kj} = \lambda^k a_{kj}, b_k = \lambda^k b_k, \; k = 1, \cdots, N, j = 1, \cdots, Z.
    \end{equation}
    }

    The algorithm of the APRFM is summarized in Algorithm~\ref{alg:aprfm-mm}.
    \begin{algorithm}[htbp]
        \caption{APRFM based on Micro-Macro Decomposition.}
        \label{alg:aprfm-mm} 
        \KwIn{The set of collocation points $\{(\bx_{\text{int}}^{k}, \bv_{\text{int}}^{k})\}_{1 \le k \le N_{\text{int}}}$, 
        $\{(\bx_{\text{bdy}}^{k}, \bv_{\text{bdy}}^{k})\}_{1 \le k \le N_{\text{bdy}}}$, 
        $\{ f_{\text{bdy}}^{k}\}_{1 \le k \le N_{\text{bdy}}}$; 
        The number of PoU functions $M^{\rho}$, $M^{g}$ and random feature functions 
        {$J^\rho, J^g$}; 
        The range of uniform distribution $[-B, B]$ for the initialization of neural network parameters; 
        The Gauss-Legendre quadrature points and weights $\{ (\omega_{j}, \bv_{j}) \}_{j=1}^{N_{q}}$. }
        \KwOut{The coefficients $\bm{\theta}^{*}$.} \BlankLine

        Initialize the weights and biases of the neural networks
        $\{\phi_{ij}^{\rho}\}, \{\phi_{ij}^{g}\}$ randomly according to the
        uniform distribution $[-B, B]$ and keep them fixed;

        Set $k_{1}= 1, k_{2}= 1$;

        \While{$k_{1}\le N_{\text{int}}$}{ \For{$i = 1, \cdots, M^{\rho}$}{ \For{$j = 1, \cdots, {J^{\rho}}$}{ Compute $^{11}s^{\text{int}, k_1}_{ij},^{21}s^{\text{int}, k_1}_{ij}$ according to \eqref{eqn:s_11} and \eqref{eqn:s_21}; } } \For{$i = 1, \cdots, M^{g}$}{ \For{$j = 1, \cdots, {J^{g}}$}{ Compute $^{12}s^{\text{int}, k_1}_{ij},^{22}s^{\text{int}, k_1}_{ij}$ according to \eqref{eqn:s_12} and \eqref{eqn:s_22}; } }

        $k_{1}= k_{1}+ 1$; }

        \While{$k_{2}\le N_{\text{bdy}}$}{ \For{$i = 1, \cdots, M^{\rho}$}{ \For{$j = 1, \cdots, {J^{\rho}}$}{ Compute $^{31}s^{\text{bdy}, k_2}_{ij}$ according to \eqref{eqn:s_31}; } } \For{$i = 1, \cdots, M^{g}$}{ \For{$j = 1, \cdots, {J^{g}}$}{ Compute $^{32}s^{\text{bdy}, k_2}_{ij}$ according to \eqref{eqn:s_32}; } }

        $k_{2}= k_{2}+ 1$; }

        Construct the matrix $\bA$ and the vector $\bm{b}$ according to~\eqref{eqn:ap-residual-matrix-A}
        and \eqref{eqn:ap-residual-vector-b};

        Solve the least squares problem~\eqref{eqn:aprfm-problem} to obtain $\bm{\theta}
        ^{*}$.

        \Return $\bm{\theta}^{*}$.
    \end{algorithm}

    \section{Numerical results}\label{sec:numerical section}

    In this section, we present the numerical results of the APRFM for the radiative transfer equation. 
    {
    To validate the effectiveness of our proposed method, we conduct all numerical experiments under the default settings of an initial weight parameter range $B = 1$, the number of Gauss-Legendre quadrature points $16$, the partition of unity function
        \begin{equation}\label{eqn:pou-fn-xd}
            \psi_{i}(\bm{y}) = \prod_{j=1}^{d}\varphi_b(\tilde{y}_{ij}),
        \end{equation}
    and the activation function $\text{tanh}(\cdot)$. 
    To enhance accuracy, we employ float64 floating-point precision and normalize all PoU functions as follows:
        \begin{equation}\label{eqn:nomalized-pou}
            \tilde{\psi_{i}}(\by) = \psi_{i}(\by) / \sum_{i=1}^{M} \psi_i(\by),
        \end{equation}
    ensuring that
   \begin{equation}\label{eqn:sum-to-one}
        \sum_{i=1}^{M} \tilde{\psi_{i}}(\bm{y}) = 1.
   \end{equation}
   For the sake of simplicity and clarity, we omit all the tildes $\tilde{}$ from $\psi_i$ throughout the paper.
    Second, we upgraded the floating-point precision from float32 to float64.}   
    The collocation points are selected as uniform grids, denoted by $(N_{x}, N_{v})$ for 1D cases and
    $(N_{x_1}, N_{x_2}, N_{v})$ for 2D cases. 
    The domain partitioning is represented as $(M_{x}, M_{v})$ for 1D cases and $(M_{x_1}, M_{x_2}, M_{v})$ for 2D cases.
    {In the numerical examples where analytical solutions are not available, 
    the reference solutions are obtained using the finite difference method.
    We compute the relative $\ell^2$ error of the solution $f(x, v)$ over all uniform grids $\{(\bx^{i}, \bv^{i})\}_{1 \le i \le I}$ with mesh size $(N_{x}, N_{v}) = (128, 256)$ for 1D cases and $(N_{x_1}, N_{x_2}, N_{v}) = (64, 64, 32)$ for 2D cases:
    \begin{equation}\label{eqn:error}
        \text{error} := \sqrt{
        \frac{\sum_{i=1}^I |f^{\text{approx}}(\bx^{i}, \bv^{i}) - f^{\text{ref}}(\bx^{i}, \bv^{i})|^2}{\sum_{i=1}^I |f^{\text{ref}}(\bx^{i}, \bv^{i})|^2}
        },
    \end{equation}
    where $f^{\text{approx}}$ denotes the solution approximated by vanilla RFM or APRFM, and $f^{\text{ref}}$ represents the reference solution.
    }

    \subsection{One-dimensional problems}
\begin{figure}[htbp]
        \centering
        \includegraphics[width=0.5\textwidth]{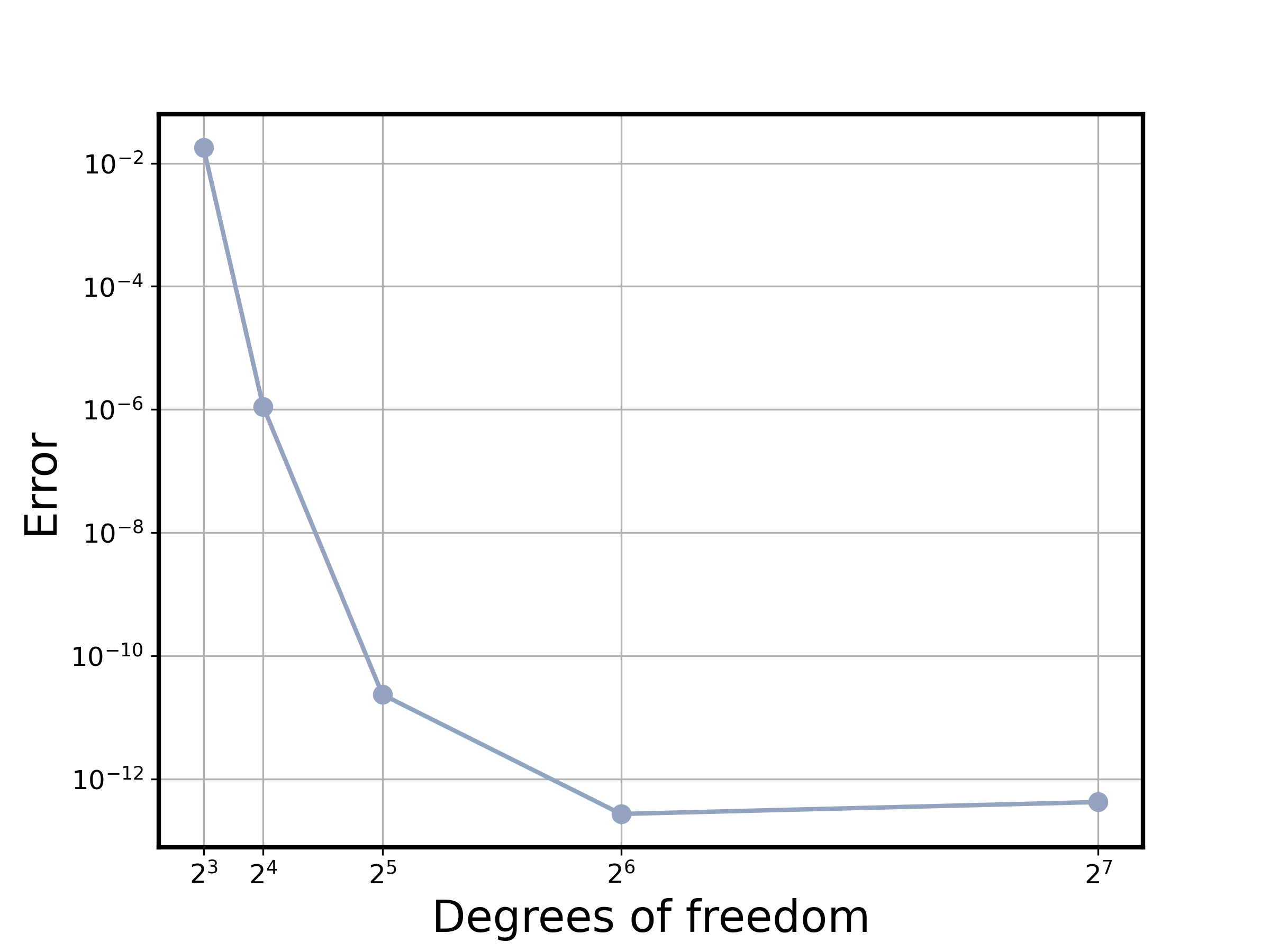}
        \caption{Relative $\ell^{2}$ error of the APRFM solution with respect to
        the degrees of freedom.}
        \label{fig:aprfm-error-semilog}
    \end{figure}

    \paragraph{Example 1}
    Let us consider the above example~\eqref{eqn:rte1d-motivation} in section~\ref{sec:motivation}:
    \begin{equation}
        \label{eqn:rte1d-ex1}
        \begin{cases}
            \displaystyle \eps v \cdot \nabla_{x}f = \avg{f}- f - \eps v, \; (x, v) \in [0, 1] \times [-1, 1], \\
            f(0, v > 0) = 1, \; f(1, v < 0) = 0,
        \end{cases}
    \end{equation}
    with the exact solution $f_{\text{ex}}= 1 - x$. 
    For the kinetic regime ($\eps=1$), 
    {we set $J^{\rho}= J^{g}= 2^{n-1}, n = 3, 4, 5, 6, 7$ and $(M_x, M_v) = (1, 1)$.} 
    Other settings are the same as Figure \ref{fig:rfm-error-semilog} in the previous section. 
    The relative $\ell^{2}$ error of the APRFM solution to the
    degrees of freedom ${Z} = 2^{n}$ is shown in Figure~\ref{fig:aprfm-error-semilog}.
    To investigate the dependency of the scale parameter in APRFM on degrees of freedom
    and collocation points, we consider $\eps = 10^{-2^i}$, $i = 1, 2, 3, 4$ with
    different random features $J^{\rho}= J^{g}= J$ and the collocation points $(N_{x}, N_{v})$. 
    We record the relative $\ell^{2}$ error of the APRFM solution with respect to the degrees of freedom 
    and the number of collocation points in Table~\ref{tab:dofs-dependency-aprfm} and Table~\ref{tab:sample-dependency-aprfm}, 
    respectively. We can observe that our APRFM is robust with respect to the
    scale parameter $\eps$ and needs fewer degrees of freedom and collocation
    points to achieve high accuracy compared with the RFM in Table~\ref{tab:dofs-dependency-rfm} and Table~\ref{tab:sample-dependency-rfm}.
    \begin{table}[tbhp]
        \caption{The dependency of the APRFM on degrees of freedom. The collocation
        points are uniform grids with $(N_{x}, N_{v}) = (128, 256)$.}
        \label{tab:dofs-dependency-aprfm}
        \centering
        \begin{tabular}{ccccccc}
            \toprule[1pt] \noalign{\smallskip} \multirow{2}*{$\eps$}                           & \multicolumn{5}{c}{$J$}  \\
                                                                                               & \multicolumn{1}{c}{$8$} & \multicolumn{1}{c}{$16$} & \multicolumn{1}{c}{$32$} & \multicolumn{1}{c}{$64$} & \multicolumn{1}{c}{$128$} \\
            \noalign{\smallskip} \midrule[1pt] \noalign{\smallskip} \multirow{1}*{{$10^{-2}$}} 
                & ${\text{7.26 e-}7}$     & ${\text{1.17 e-}12}$      & ${\text{5.71 e-}14}$      & ${\text{2.69 e-}14}$      & ${\text{4.72 e-}14}$       \\
            \multirow{1}*{{$10^{-4}$}}                                                         
                & ${\text{7.27 e-}7}$     & ${\text{1.55 e-}12}$      & ${\text{9.76 e-}14}$      & ${\text{5.75 e-}14}$      & ${\text{4.70 e-}14}$       \\
            \multirow{1}*{{$10^{-8}$}}                                                         
                & ${\text{7.45 e-}7}$     & ${\text{1.02 e-}12}$      & ${\text{5.57 e-}14}$      & ${\text{4.28 e-}14}$      & ${\text{4.35 e-}14}$       \\
            \multirow{1}*{{$10^{-16}$}}                                                        
                & ${\text{7.54 e-}7}$     & ${\text{2.67 e-}12}$      & ${\text{1.55 e-}14}$      & ${\text{1.66 e-}14}$      & ${\text{3.11 e-}14}$       \\
            \noalign{\smallskip} \bottomrule[1pt]
        \end{tabular}
    \end{table}
    \begin{table}[tbhp]
        \caption{The dependency of the APRFM on the number of collocation points.
        The number of random features $J = 128$.}
        \label{tab:sample-dependency-aprfm}
        \centering
        \begin{tabular}{cccccc}
            \toprule[1pt] \noalign{\smallskip} \multirow{2}*{$\eps$}                           & \multicolumn{4}{c}{$(N_{x}, \, N_{v})$} \\
                                                                                               & \multicolumn{1}{c}{$(16, \, 32)$}      & \multicolumn{1}{c}{$(32, \, 64)$} & \multicolumn{1}{c}{$(64, \, 128)$} & \multicolumn{1}{c}{$(128, \, 256)$} \\
            \noalign{\smallskip} \midrule[1pt] \noalign{\smallskip} \multirow{1}*{{$10^{-2}$}} 
             & ${\text{4.38 e-}10}$               & ${\text{2.05 e-}11}$                       & ${\text{5.99 e-}13}$                   & ${\text{4.72 e-}14}$         \\
            \multirow{1}*{{$10^{-4}$}}                                                         
             & ${\text{4.34 e-}10}$               & ${\text{2.16 e-}11}$                       & ${\text{6.21 e-}13}$                   & ${\text{4.70 e-}14}$         \\
            \multirow{1}*{{$10^{-8}$}}                                                         
             & ${\text{2.37 e-}8}$                & ${\text{2.14 e-}11}$                       & ${\text{6.37 e-}13}$                   & ${\text{4.35 e-}14}$         \\
            \multirow{1}*{{$10^{-16}$}}                                                        
             & ${\text{2.31 e-}8}$                    & ${\text{2.18 e-}11}$                   & ${\text{6.01 e-}13}$                   & ${\text{3.11 e-}14}$          \\
            \noalign{\smallskip} \bottomrule[1pt]
        \end{tabular}
    \end{table}

\begin{figure}[htbp!]
        \centering
        \subfigure[Reference solution v.s. APRFM solution ($\eps = 1$).]{ 
            \includegraphics[width=\textwidth]{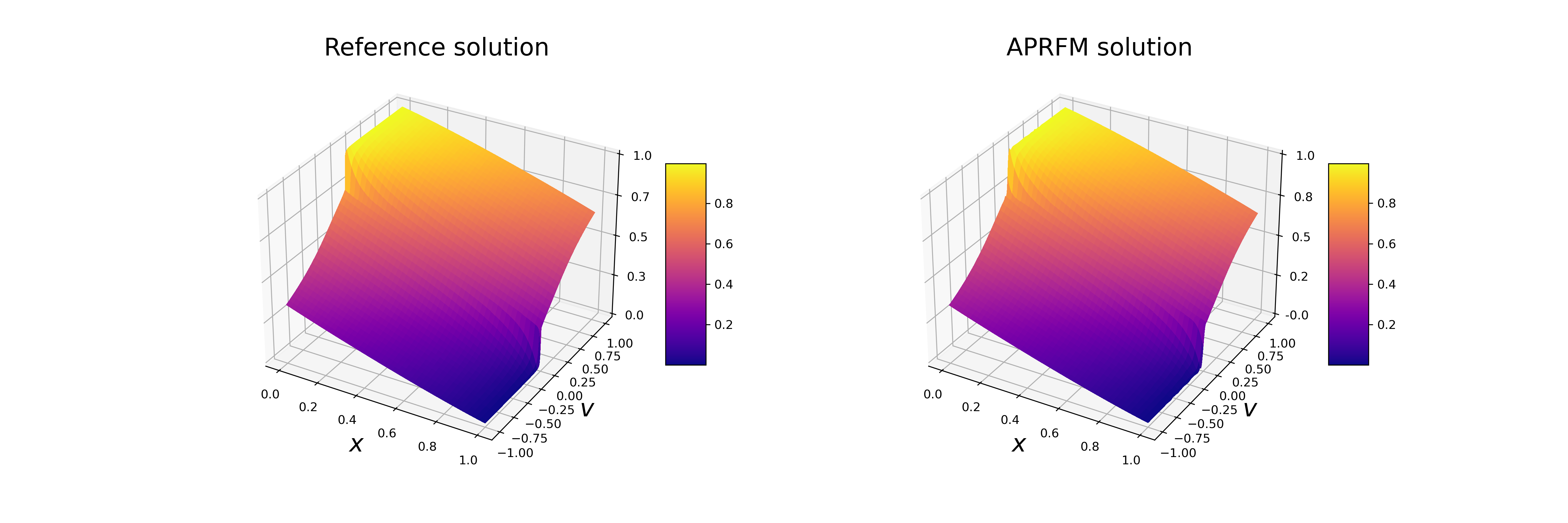}
        } \subfigure[Reference solution v.s. APRFM solution ($\eps = 5 \times 10^{-1}$).]{ 
            \includegraphics[width=\textwidth]{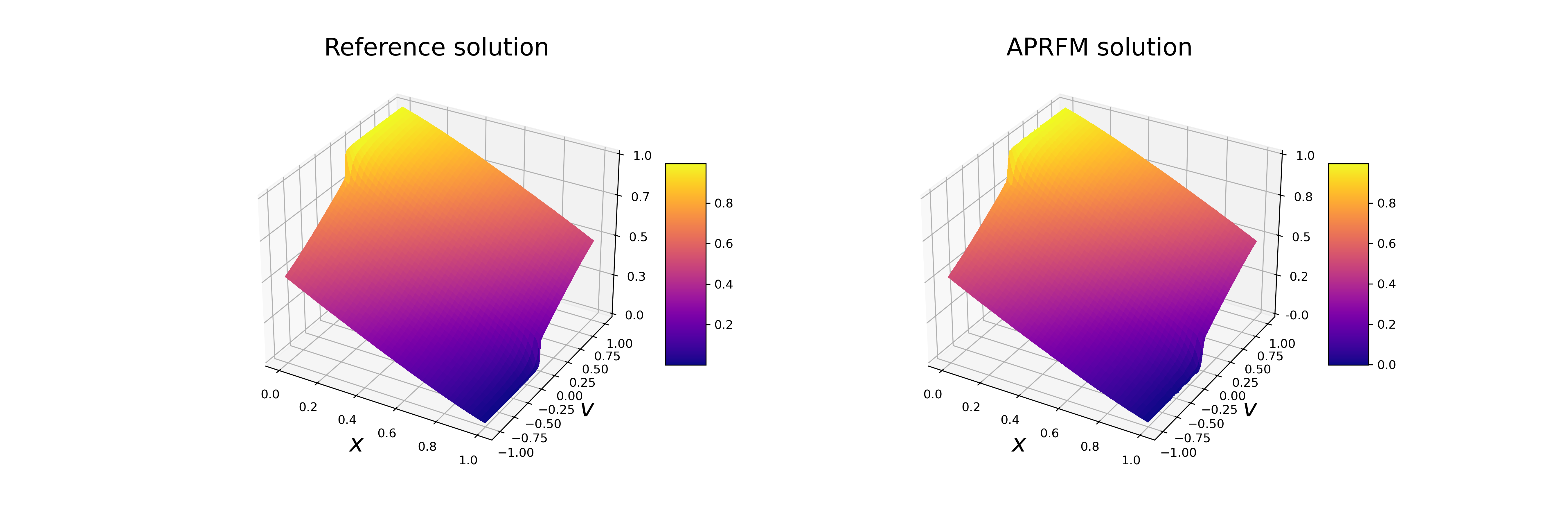}
        }
        \caption{Plots of $f(x, v)$ in $(x, v) \in [0, 1] \times [-1, 1]$.}
        \label{fig:aprfm1d-ex2}
    \end{figure}

    \paragraph{Example 2}
    In this example, we consider the kinetic regime ($\eps = 1$) and
    intermediate regime ($\eps = 5 \times 10^{-1}$) in 1D slab geometry without the presence of a source term:
    \begin{equation}
        \label{eqn:rte1d-ex2}
        \begin{cases}
            \displaystyle \eps v \cdot \nabla_{x}f = \avg{f}- f, \; (x, v) \in [0, 1] \times [-1, 1], \\
            f(0, v > 0) = 1, \; f(1, v < 0) = 0.
        \end{cases}
    \end{equation}
    We plot the solution obtained by our APRFM and reference solution in Figure~\ref{fig:aprfm1d-ex2}.
    {We set $J^{\rho}= 64, J^{g}= 128$.
    The number of interior collocation points is $(N_{x}, N_{v}) = (128, 256)$.
    Additionally, the domain is partitioned as $(M_{x}, M_{v}) = (2, 4)$, i.e., $M^{\rho}= 2, M^{g}= 8$.
    The relative $\ell^{2}$ error of our APRFM is $7.08 \times 10^{-3}$ for $\eps = 1$ and $5.56 \times 10^{-3}$ for $\eps = 5 \times 10^{-1}$, respectively. }
    It can be observed that the errors of the APRFM are predominantly concentrated at $x = 0, 1$ and $v = 0$.
    
\begin{figure}[htbp!]
        \centering
        \includegraphics[width=0.8\textwidth]{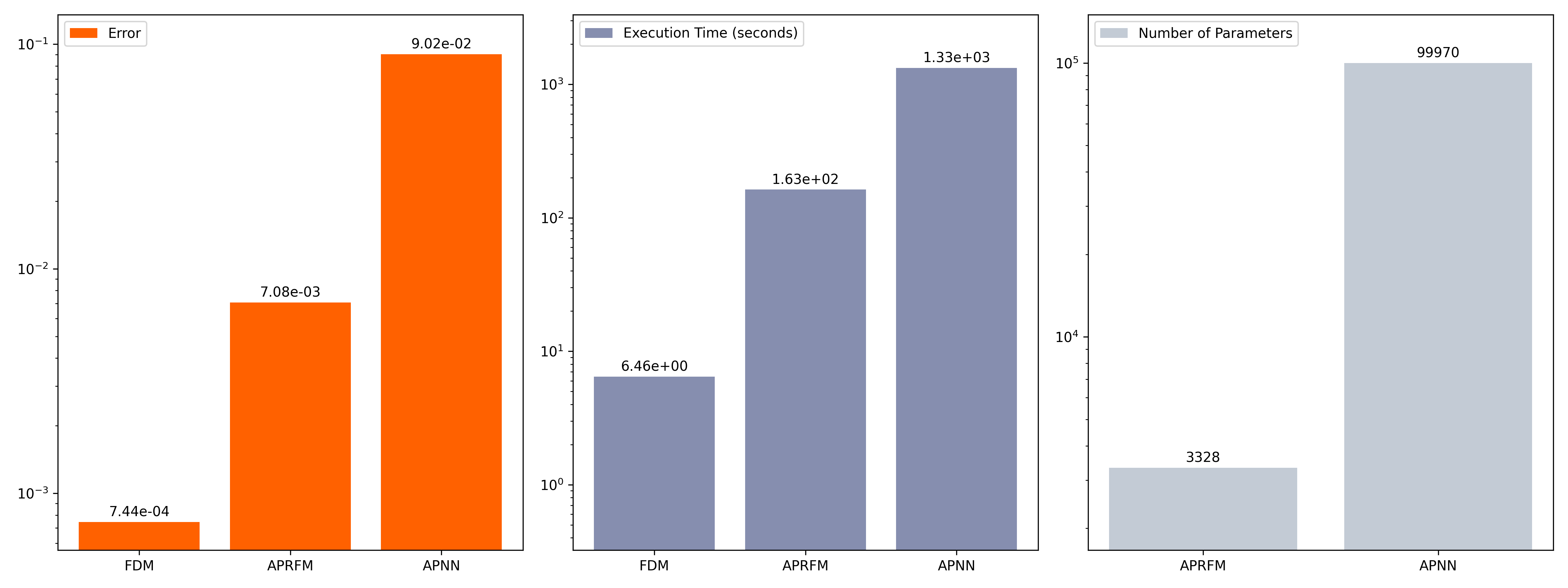}
        \caption{{The comparison of APRFM and APNN in 1D case.}}
        \label{fig:comparison}
    \end{figure}

    In the kinetic case ($\eps = 1$), we also compared our proposed APRFM with the {finite difference method (FDM) and} APNN method~\cite{wuAPNN,lu2022solving}. 
    As shown in Figure~\ref{fig:comparison}, our approach not only achieves lower relative $\ell^{2}$ error compared to the APNN method, but more importantly, it significantly outperforms in terms of computational efficiency, being an order of magnitude faster and requiring considerably fewer parameters. 
    Owing to the uniform accuracy of our APRFM, the comparative results hold similarly for cases where $\eps \ll 1$.
    {For this 1D example, our method does not provide a significant advantage in terms of computation time compared to the traditional finite difference method. 
    In this example, due to the large number of collocation points and the extensive partitioning of the domain, the assembly of the matrix took 144.4 seconds, accounting for $88\%$ of the total computation time.
    However, it is important to highlight that our APRFM demonstrates notable improvements in computation time for 2D cases, as detailed in Example 4.}
    This comparative experiment was conducted on a machine equipped with an NVIDIA A800 80G GPU and a 48-core Intel Xeon Gold 6342 CPU.
    The neural networks for $\rho$ and $g$ in the APNN method are 4-layer fully connected architectures, each with 128 units per hidden layer, 
    utilizing the $\tanh$ activation function. 
    We trained the APNN method for 50000 iterations using the Adam optimizer, 
    with a batch size of 8196 for interior points and 2048 for boundary points during the training process. 
    {The number of inner-layer network parameters in our APRFM that do not participate in training is 1152. Other settings are the same as those in Example 2.}

\begin{figure}[htbp!]
        \centering
        \includegraphics[width=0.5\textwidth]{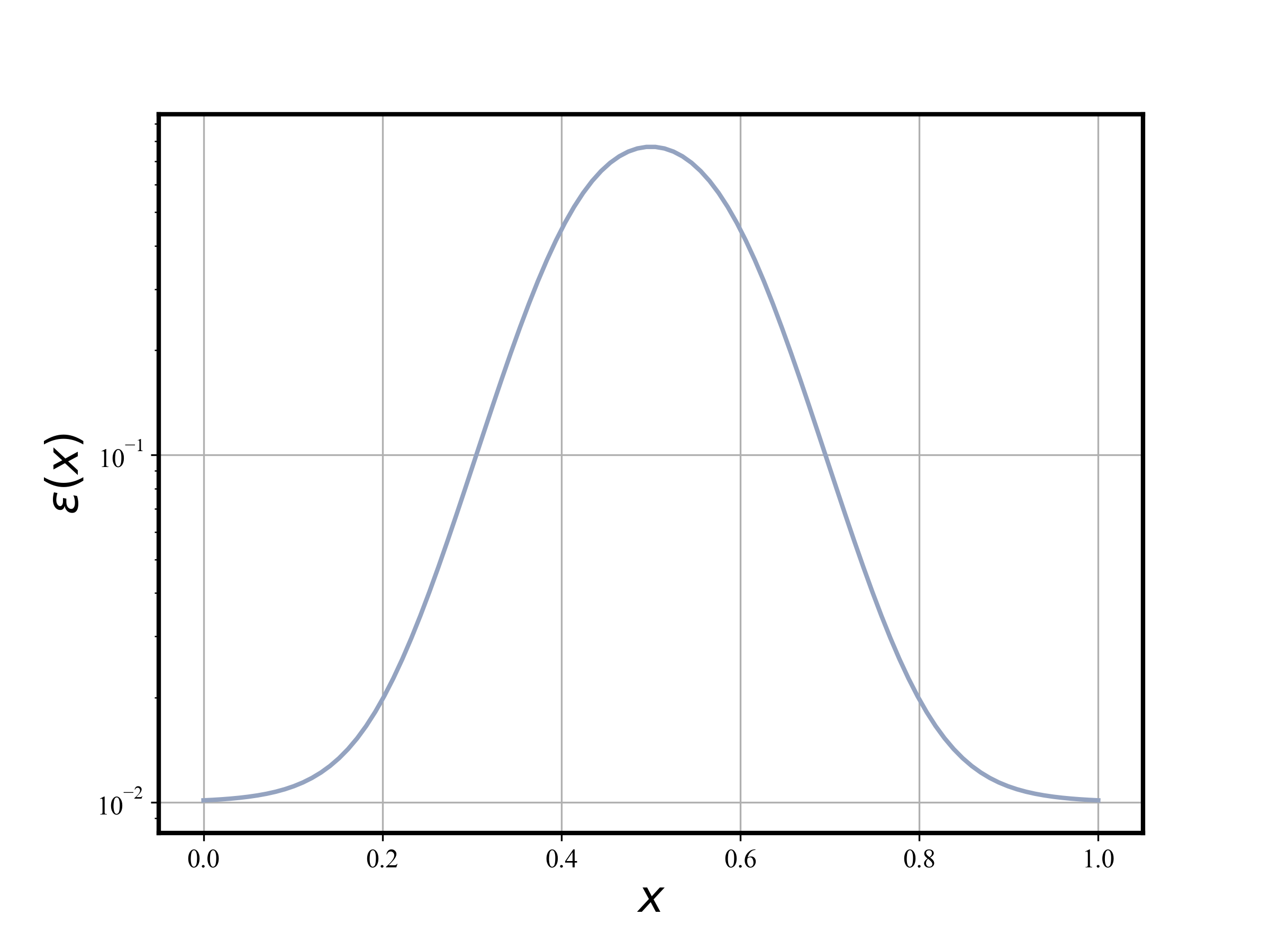}
        \caption{Plot of $\eps(x)$.}
        \label{fig:mixed-scale}
    \end{figure}

    \paragraph{Example 3}
    The following example involves a mixed-scale RTE:
    \begin{equation}
        \label{eqn:rte1d-ex3}
        \begin{cases}
            {}\displaystyle v \cdot \nabla_{x}f = \frac{1}{\eps(x)}\left ( \avg{f}- f \right ), \; (x, v) \in [0, 1] \times [-1, 1], \\
            f(0, v > 0) = 0.5, \; f(1, v < 0) = 0.
        \end{cases}
    \end{equation}
    This case can verify the adaptability of our APRFM under mixed multi-scale conditions. 
    Here, $\eps(x)$ is a function that depends on the spatial
    variable and smoothly transitions from {$\mathcal{O}(10^{-2})$} to $\mathcal{O}(1)$, 
    defined as follows:
    \begin{equation}
        \eps(x) = 10^{-2}+ \frac{1}{2}\left [ \tanh (6.5 -11 x) + \tanh (11 x - 4
        .5) \right ].
    \end{equation}
    The plot of the function $\eps(x)$ is presented in Figure~\ref{fig:mixed-scale}.
    
    {
    For this case, we briefly derive the corresponding micro-macro decomposition system.
    First, we decompose $f(x, v)$ as 
    \begin{equation}
    \label{eqn:mixed-decompose}
        f(x, v) = \rho(x) + \eps(x) g(x, v), 
    \end{equation}
    with $\avg{g} = 0$.
    Substituting this into equation~\eqref{eqn:rte1d-ex3}, we rewrite it as
    \begin{equation}
    \label{eqn:mixed-formula}
        v \cdot \nabla_{x} (\rho(x) + \eps(x) g(x, v)) + g(x, v) = 0.
    \end{equation}
    Applying the orthogonal projection operators $\Pi: \Pi(\cdot)(v) = \avg{\cdot}$ and $\text{Id} - \Pi$ to equation~\eqref{eqn:mixed-formula}, we derive the micro-macro system for this mixed-scale problem:
    \begin{equation}
        \label{eqn:mixed-micro-macro}
        \begin{cases}
            \avg{v \cdot \nabla_{x} (\eps(x) g)} = 0, \\
            v \cdot \nabla_{x} \rho + (\text{Id}- \Pi) (v \cdot \nabla_{x} (\eps(x) g)) + g = 0. 
        \end{cases}
    \end{equation}
    Thus, we can employ our APRFM to solve the micro-macro system~\eqref{eqn:mixed-micro-macro}.}  
    
\begin{figure}[htbp!]
        \centering
        \includegraphics[width=\textwidth]{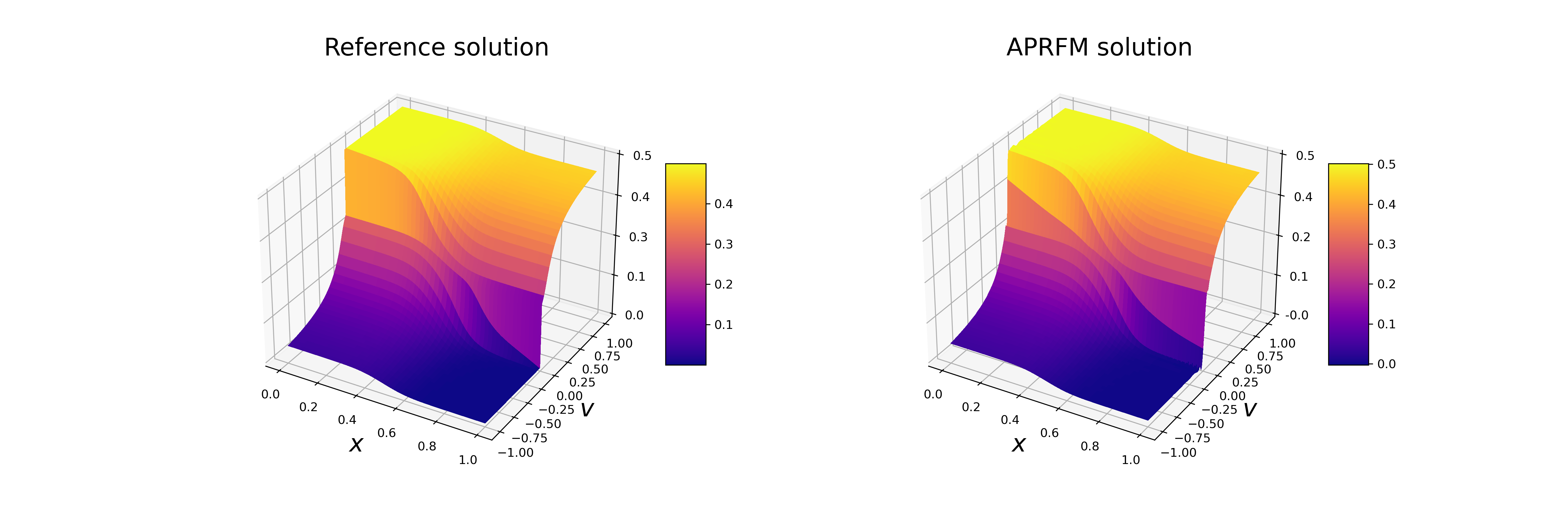}
        \caption{Reference solution v.s. APRFM solution.}
        \label{fig:aprfm1d-ex3-sigma_s}
    \end{figure}

    We plot the solution obtained by our APRFM and reference solution in Figure~\ref{fig:aprfm1d-ex3-sigma_s}.
    {We set $J^{\rho}= 64, J^{g}= 128$.
    The number of collocation points is $(N_{x}, N_{v}) = (128, 256)$.
    The domain is partitioned as $(M_{x}, M_{v}) = (2, 4)$, i.e., $M^{\rho}= 2, M^{g}= 8$.  
    The relative $\ell^{2}$ error of our APRFM is $1.72 \times 10^{-2}$.}

    \subsection{Two-dimensional problems}
\begin{figure}[htbp!]
        \centering
        \subfigure[Reference solution v.s. APRFM solution ($\eps = 1$).]{
        \includegraphics[width=\textwidth]{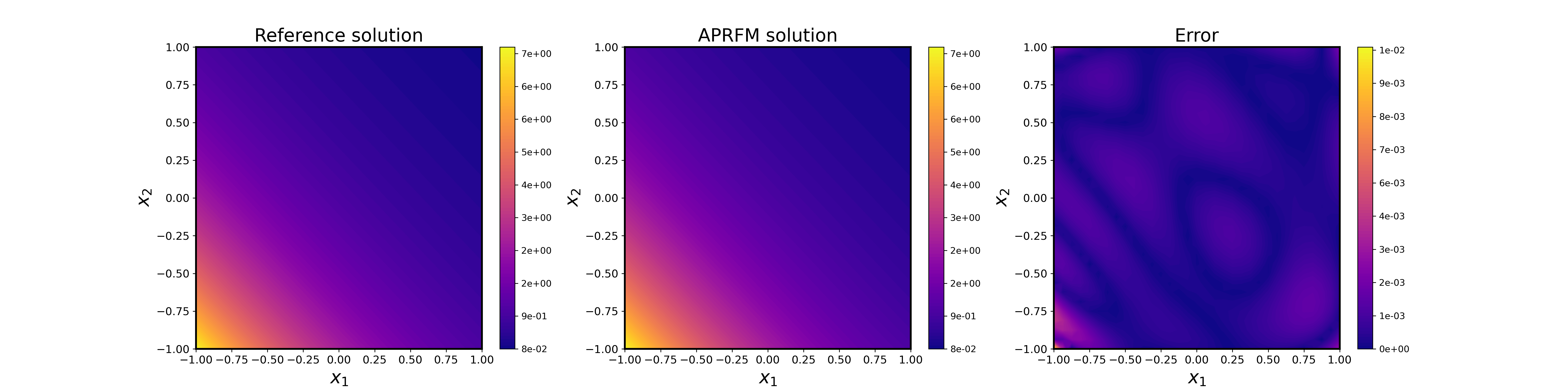} } 
        \subfigure[Reference solution v.s. APRFM solution ($\eps = 10^{-3}$).]{
        \includegraphics[width=\textwidth]{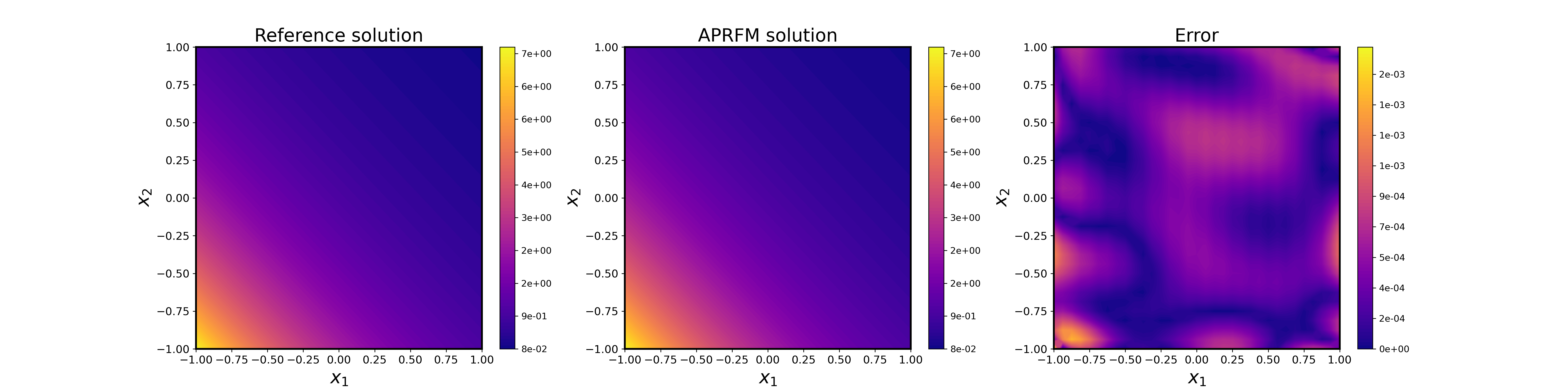} }
        \caption{Plots of $\rho(x_1, x_2)$ in $(x_1, x_2) \in [-1, 1] \times [-1, 1]$.}
        \label{fig:aprfm2d-ex4}
    \end{figure}

    \paragraph{Example 4}
    In this subsection, we consider a 2D square domain 
    with $\bx = (x_{1}, x_{2}) \in [-1, 1]^{2}, \bv = (\cos(\alpha), \sin(\alpha))$.
    \begin{equation}
        \label{eqn:rte2d-ex1}
        \begin{cases}
            {}\displaystyle \eps \bv \cdot \nabla_{\bx}f(\bx, \bv) = \frac{1}{2 \pi}\int_{|\bv| = 1}f(\bx, \bv') \diff{\bv'}- f + \eps^{2}G(\bx, \bv),  \\
            f(-1, x_{2}, \alpha) = e^{1 - x_2}, \; \alpha \in [0, \pi/2] \cup [3\pi/2, 2\pi],                                                                                                    \\
            f(1, x_{2}, \alpha) = e^{- 1 - x_2}, \; \alpha \in [\pi/2, 3\pi/2],                                                                                                                  \\
            f(x_{1}, -1, \alpha) = e^{- x_1 + 1}, \; \alpha \in [0, \pi],                                                                                                                        \\
            f(x_{1}, 1, \alpha) = e^{- x_1 - 1}, \; \alpha \in [\pi, 2\pi].
        \end{cases}
    \end{equation}
    {Here, the source function $G(\bx, \bv) = (-\cos(\alpha) - \sin(\alpha)) \exp(- x_1 - x_2) / \eps$.
    In this example, the exact solution is given by $f_{\text{ex}} = \exp(- x_1 - x_2)$.
    We plot the density function $\rho(\bx)$ obtained by our APRFM and reference solution 
    in Figure~\ref{fig:aprfm2d-ex4}.  
    We set $J^{\rho}= 32, J^{g}= 32$ and the number of collocation points is $(N_{x_1}, N_{x_2}, N_{v}) = (32, 32, 64)$.
    The domain is partitioned as $(M_{x_1}, M_{x_2}, M_{v}) = (1, 1, 1)$, i.e., $M^{\rho}= 1, M^{g}= 1$.
    The relative $\ell^{2}$ error of density function $\rho(\bx)$ obtained by 
    our APRFM is $3.48 \times 10^{-4} (\eps = 1)$ and $2.28 \times 10^{-4} (\eps = 10^{-3})$.
    For $\eps = 1$ of this example, the computation time for our APRFM was 6.29 seconds, while the computation time for the finite difference method was 111.20 seconds.
    It should be noted that as the number of collocation points, partition regions, and random feature functions increase, the computation time for the APRFM will correspondingly increase.
    }
    
\begin{figure}[htbp!]
        \centering
        \subfigure[Reference solution v.s. APRFM solution ($\eps = 1$).]{
        \includegraphics[width=\textwidth]{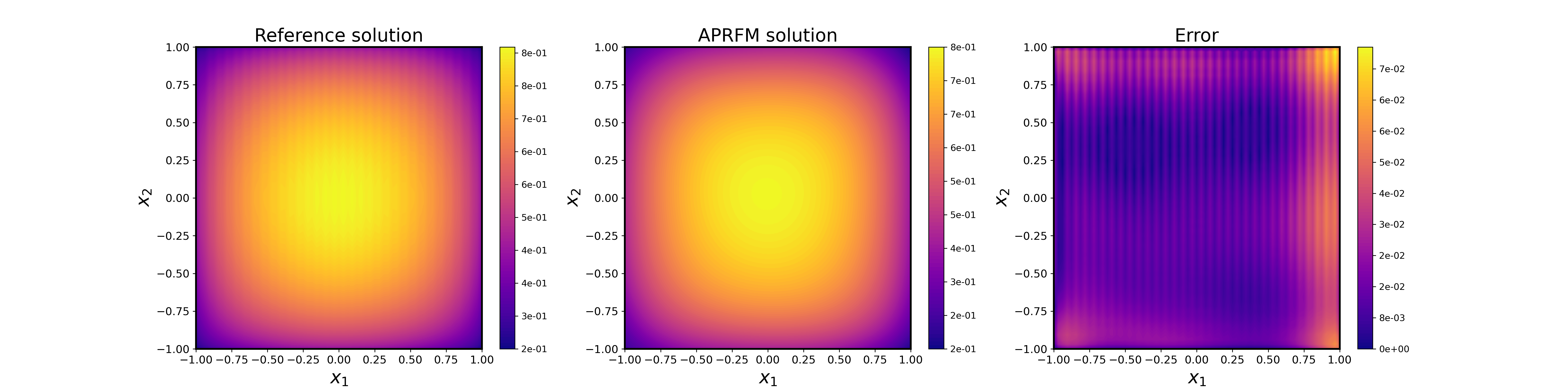} } \subfigure[Reference
        solution v.s. APRFM solution ($\eps = 10^{-1}$).]{

        \includegraphics[width=\textwidth]{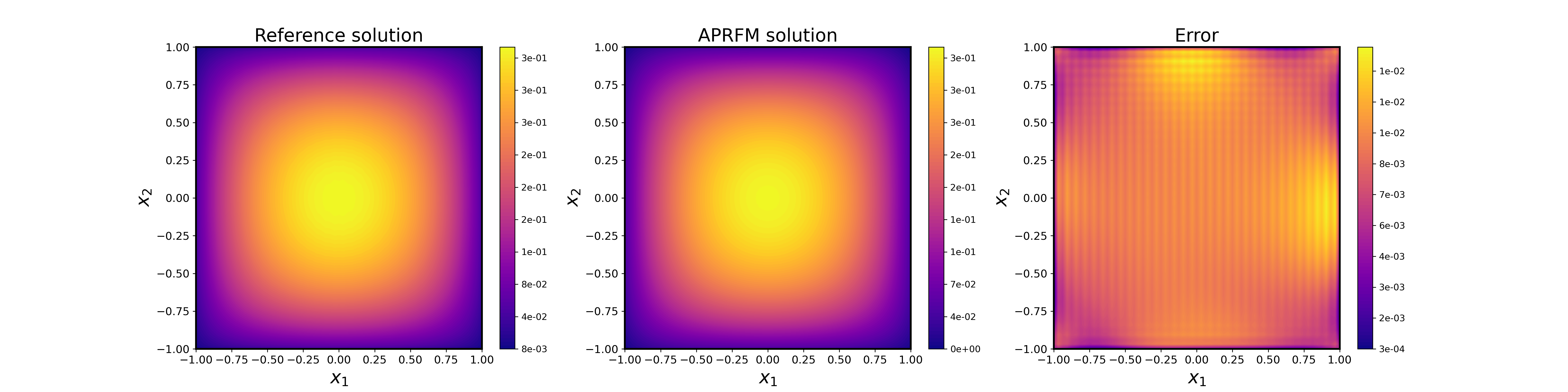} }
        \caption{Plots of $\rho(x_1, x_2)$ in $(x_1, x_2) \in [-1, 1] \times [-1, 1]$ by
        APRFM and finite difference method.}
        \label{fig:aprfm2d-ex5}
    \end{figure}

    \paragraph{Example 5}
    Next, we consider a square domain with $\bx = (x_{1}, x_{2}) \in [-1, 1]^{2},
    \bv = (\cos(\alpha), \sin(\alpha))$.
    \begin{equation}
        \label{eqn:rte2d-ex2}
        \begin{cases}
            \displaystyle \eps \bv \cdot \nabla_{\bx}f(\bx, \bv) = \frac{1}{2 \pi}\int_{|\bv| = 1}f(\bx, \bv') \diff{\bv'}- f + \eps^{2}G(\bx, \bv),  \\
            f(-1, x_{2}, \alpha) = 0, \; \alpha \in [0, \pi/2] \cup [3\pi/2, 2\pi],                                                                   \\
            f(1, x_{2}, \alpha) = 0, \; \alpha \in [\pi/2, 3\pi/2],                                                                                   \\
            f(x_{1}, -1, \alpha) = 0, \; \alpha \in [0, \pi],                                                                                         \\
            f(x_{1}, 1, \alpha) = 0, \; \alpha \in [\pi, 2\pi].
        \end{cases}
    \end{equation}
    {
    Here, the source function $G(\bx, \bv) = 1/2$.
    We plot the density function $\rho(\bx)$ obtained by our APRFM and reference solution 
    in Figure~\ref{fig:aprfm2d-ex5}.
    We set $J^{\rho}= 64, J^{g}= 128$ and the number of interior collocation points is $(N_{x_1}, N_{x_2}, N_{v}) = (32, 32, 32)$. 
    The domain is partitioned as $(M_{x_1}, M_{x_2}, M_{v}) = (1, 1, 4)$, resulting in $M^{\rho}= 1, M^{g}= 4$.
    The relative $\ell^{2}$ error of our APRFM is $3.42 \times 10^{-2} (\eps = 1)$ and $4.43 \times 10^{-2} (\eps = 10^{-1})$. 
    Additionally, we have documented the relative $\ell^{2}$ errors corresponding to different PoU configurations in Table~\ref{tab:pou-dependency-rfm2d}.
    In this case, it has been observed that simply increasing the number of partitions does not necessarily lead to better results.
    }
    
    \begin{table}[tbhp]
        \caption{The dependency of the RFM on PoU. The collocation
        points are uniform grids with $(N_{x_1}, N_{x_2}, N_{v}) = (32, 32, 32)$ and the number of random feature functions $J^{\rho}= 64, J^{g}= 128$.}
        \label{tab:pou-dependency-rfm2d}
        \centering
        \begin{tabular}{cccccc}
            \toprule[1pt] \noalign{\smallskip} \multirow{2}*{$\eps$}                           & \multicolumn{4}{c}{$(M_{x_1}, \, M_{x_2}, \, M_{v})$} \\
                                                                                               & \multicolumn{1}{c}{$(1, \, 1, \, 1)$}   & \multicolumn{1}{c}{$(1, \, 1, \, 2)$} & \multicolumn{1}{c}{$(1, \, 1, \, 4)$} & \multicolumn{1}{c}{$(1, \, 1, \, 8)$} \\
            \noalign{\smallskip} \midrule[1pt] \noalign{\smallskip} \multirow{1}*{{$1$}} 
            & ${\text{1.40 e-}1}$        & ${\text{5.82 e-}2}$      & ${\text{3.42 e-}2}$      & ${\text{4.80 e-}2}$       \\
            \multirow{1}*{{$10^{-1}$}}   
            & ${\text{4.23 e-}2}$        & ${\text{4.83 e-}2}$      & ${\text{4.43 e-}2}$      & ${\text{1.56 e-}1}$       \\    
            \noalign{\smallskip} \bottomrule[1pt]
        \end{tabular}
    \end{table}

\begin{figure}[htbp!]
        \centering
        \subfigure[Reference solution v.s. APRFM solution ($\eps = 1$).]{
        \includegraphics[width=\textwidth]{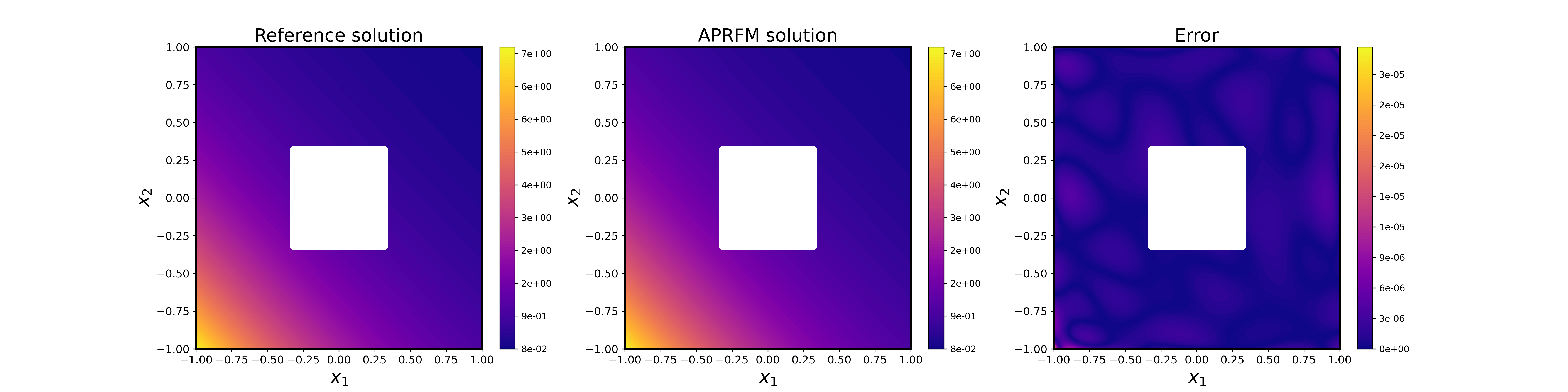} }
        \subfigure[Reference solution v.s. APRFM solution ($\eps = 5 \times 10^{-3}$).]{
        \includegraphics[width=\textwidth]{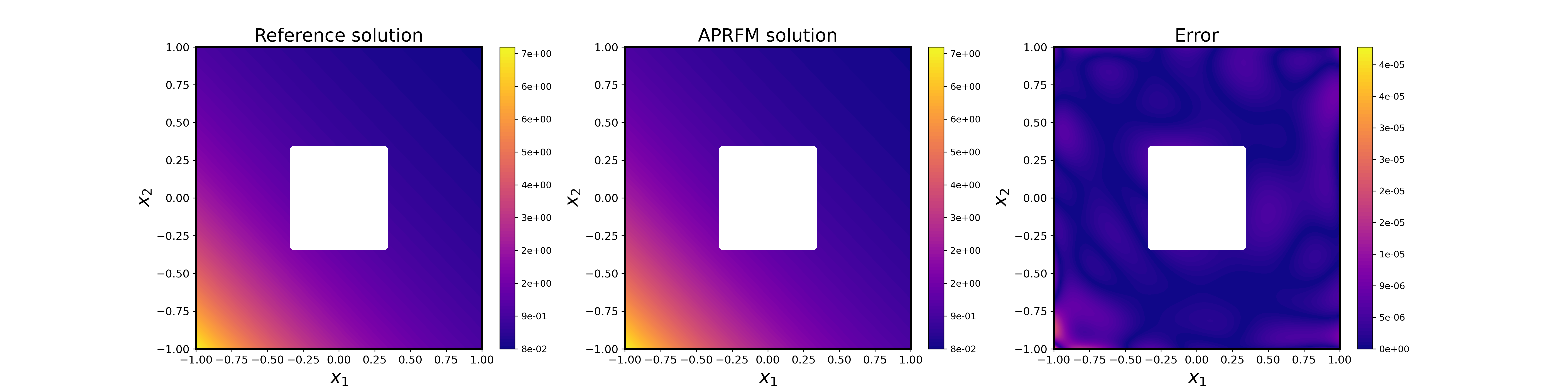} }
        \caption{Plots of $\rho(x, y)$ in $(x_1, x_2) \in [-1, 1] \times [-1, 1]$.}
        \label{fig:aprfm2d-ex6}
    \end{figure}

\begin{figure}[htbp!]
        \centering
        \subfigure[Reference solution v.s. APRFM solution ($\eps = 1$).]{
        \includegraphics[width=\textwidth]{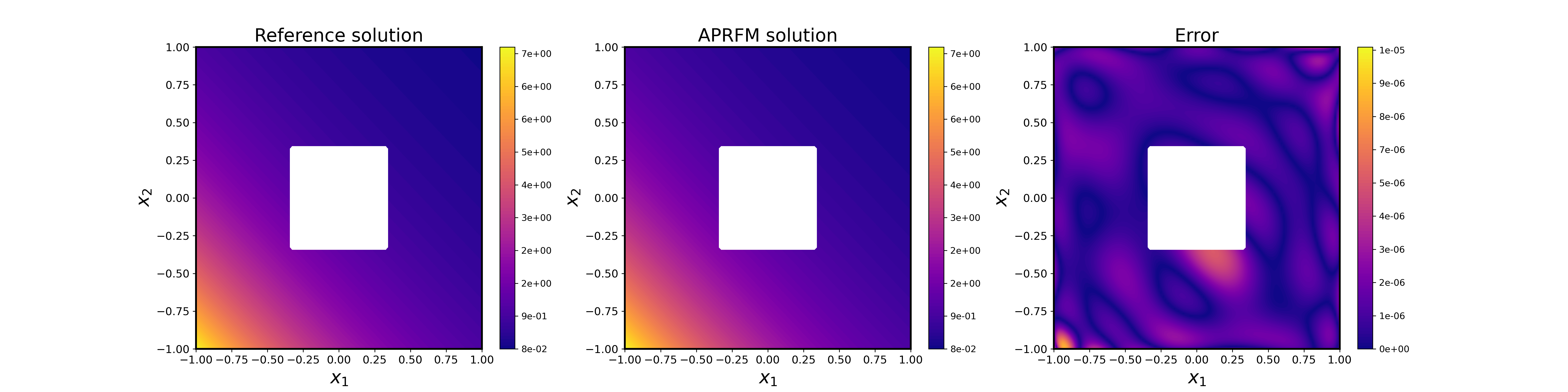} }
        \subfigure[Reference solution v.s. APRFM solution ($\eps = 5 \times 10^{-3}$).]{
        \includegraphics[width=\textwidth]{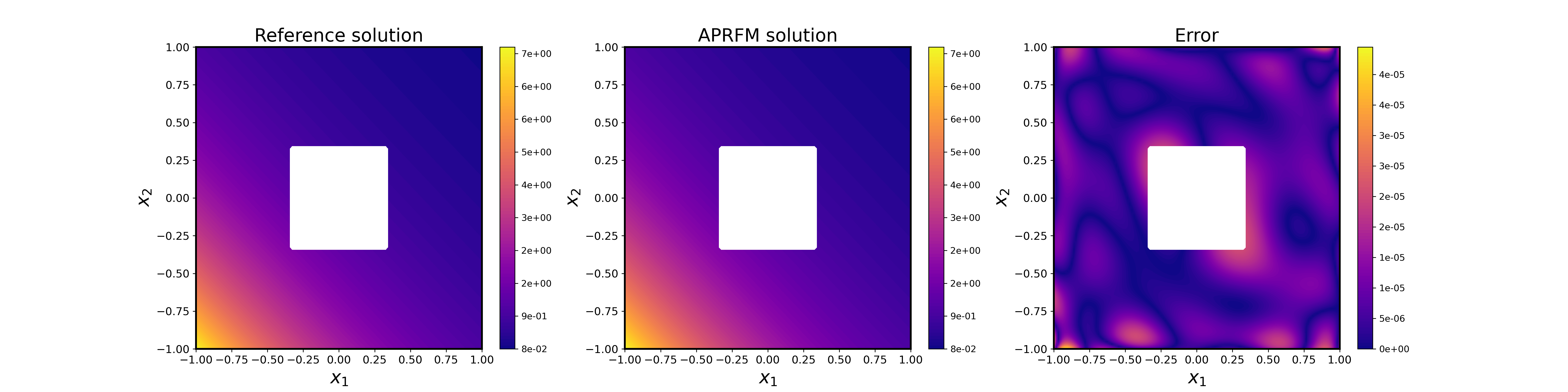} }
        \caption{Plots of $\rho(x, y)$ in $(x_1, x_2) \in [-1, 1] \times [-1, 1]$.}
        \label{fig:aprfm2d-ex6-sin}
    \end{figure}

    \paragraph{Example 6}
    Finally, we focus on a hollow annular domain with 
    $\bx = (x_{1}, x_{2}) \in [-1, 1]^{2} \setminus (-1/3, 1/3)^{2}, 
    \bv = (\cos(\alpha), \sin(\alpha))$.
    \begin{equation}
        \label{eqn:rte2d-ex3}
        \begin{cases}
            \displaystyle \eps \bv \cdot \nabla_{\bx}f(\bx, \bv) = \frac{1}{2 \pi}\int_{|\bv| = 1}f(\bx, \bv') \diff{\bv'}- f + \eps^{2}G(\bx, \bv),  \\
            f(-1, x_{2}, \alpha) = e^{1 - x_2}, \, f(\frac{1}{3}, x_{2}, \alpha) = e^{-\frac{1}{3} - x_2}, \; \alpha \in [0, \pi/2] \cup [3\pi/2, 2\pi],                                         \\
            f(1, x_{2}, \alpha) = e^{- 1 - x_2}, \, f(-\frac{1}{3}, x_{2}, \alpha) = e^{\frac{1}{3} - x_2}, \; \alpha \in [\pi/2, 3\pi/2],                                                       \\
            f(x_{1}, -1, \alpha) = e^{- x_1 + 1}, \, f(x_{1}, \frac{1}{3}, \alpha) = e^{-x_1 - \frac{1}{3}}, \; \alpha \in [0, \pi],                                                             \\
            f(x_{1}, 1, \alpha) = e^{- x_1 - 1}, \, f(x_{1}, -\frac{1}{3}) = e^{-x_1 + \frac{1}{3}}, \; \alpha \in [\pi, 2\pi] .
        \end{cases}
    \end{equation}
    {
    Here, the source function $G(\bx, \bv) = (-\cos(\alpha) - \sin(\alpha)) \exp(- x_1 - x_2) / \eps$.
    As before, the exact solution is given by {$f_{\text{ex}} = \exp(- x_1 - x_2)$}.
    We plot the density function $\rho(\bx)$ obtained by our APRFM and reference solution 
    in Figure~\ref{fig:aprfm2d-ex6}.
    We set $J^{\rho}= 64, J^{g}= 128$ and the number of collocation points is $(N_{x_1}, N_{x_2}, N_{v}) = (32, 32, 64)$.
    The domain is partitioned as $(M_{x_1}, M_{x_2}, M_{v}) = (1, 1, 4)$, resulting in $M^{\rho}= 1, M^{g}= 4$.
    The relative $\ell^{2}$ error of our APRFM is $7.54 \times 10^{-7} (\eps = 1)$ and 
    $1.60 \times 10^{-6} (\eps = 5 \times 10^{-3})$.
    Besides, we replaced the activation function with $\sin (\pi x)$, while keeping all other settings unchanged. 
    The corresponding density function $\rho$ is shown in Figure~\ref{fig:aprfm2d-ex6-sin}, 
    with the relative $\ell^{2}$ error of $7.34 \times 10^{-7} (\eps = 1)$ and 
    $4.50 \times 10^{-6} (\eps = 5 \times 10^{-3})$.
    }

    \section{Conclusions}   
    In this paper, we proposed a novel efficient Random Feature Method based on 
    micro-macro decomposition for efficiently solving multiscale radiative transfer equations. Motivated by a simple example, we observed that the vanilla RFM is unsuitable for the kinetic regime,  while the proposed APRFM effectively addresses the RTE in this context. Extensive numerical experiments have been conducted to verify the effectiveness of the APRFM for both 1D and 2D RTEs.
    The numerical results of the APRFM demonstrate that it can achieve high accuracy 
    with fewer degrees of freedom and collocation points compared with the vanilla RFM.
    The APRFM is robust with respect to the scale parameter $\eps$ and can be applied to mixed multi-scale RTEs. The APRFM is a promising method for solving multiscale radiative transfer equations, and it has the potential to be extended to time-dependent and nonlinear kinetic equations, as well as complex geometries. The least squares problem in the APRFM, however, is more challenging to solve compared to traditional numerical methods. This involves factors such as the random initialization of network parameters, the selection of activation functions, and the choice of the PoU function, efficient sampling strategies, and domain decomposition techniques. Besides, the analysis of the convergence and stability of the APRFM is still an open question. In future work, we plan to enhance the efficiency and accuracy of the APRFM in solving kinetic equations and to extend the method to time-dependent and nonlinear kinetic equations. We also aim to investigate the convergence and stability of the APRFM and to explore the application of the APRFM to other multiscale problems.

    \section*{Acknowledgments}
    Jingrun Chen is partially supported by the NSFC Major Research Plan - Interpretable and General Purpose Next-generation Artificial Intelligence (Nos. 92270001 and 92370205), NSFC grant 12425113, and Key Laboratory of the Ministry of Education for Mathematical Foundations and Applications of Digital Technology, University of Science and Technology of China.
    Zheng Ma is supported by NSFC No. 12201401, No. 92270120 and Beijing Institue of Applied Physics and Computational Mathematics funding HX02023-60.

    \bibliographystyle{unsrt} 
    \bibliography{references}

\end{document}